\documentclass{amsart}
\usepackage{amsfonts}
\usepackage{latexsym}
\usepackage{amssymb}
\usepackage{amsmath}
\usepackage{todonotes}
\usepackage{color}
\usepackage{pdfpages}
\usepackage{graphicx}
\usepackage{subcaption}

\def\N{\mathbb{N}}

\begin{document}


\title{El Logaritmo Integral: Origen y aplicaciones en las matemáticas el siglo XIX.}
\title{El Logaritmo Integral: Números primos y algo más.}

\author[M.A. Crespo]{Miguel Angel Crespo Mir}
\address{Investigador independiente. Colaborador del IUMA.}
\email[(Miguel Angel Crespo Mir)]{mcrespomir@icloud.com}
\thanks{Partially supported by IUMA}
\date{\today}

\author[J.\,Bernu\'es]{Julio Bernu\'es}
\address{Instituto Universitario de Matemáticas y Aplicaciones (IUMA), Universidad de Zaragoza, 50009 Zaragoza (Spain)}
\email[(Julio Bernu\'es)]{bernues@unizar.es}

\begin{abstract} En este artículo mostramos la relevancia que la función logaritmo integral tuvo en el desarrollo de las matemáticas de la primera mitad del siglo XIX. Su importancia implicará a matemáticos de primer nivel como Euler, Gauss, Bessel, Riemann.

La visión que presentamos es fruto del estudio pormenorizado de las fuentes originales 
que nos permite establecer la línea temporal de cómo se sucedieron los avances. En particular, nuestro estudio reivindica las aportaciones de Bessel (en colaboración con Gauss) que el paso del tiempo parecía haber diluido. 

El logaritmo integral en la mente de Gauss jugará un papel fundamental que le hará anticiparse a su tiempo en campos tan diversos como el análisis complejo (teorema de Cauchy),  métodos numéricos de cuadratura (cuadratura de Gauss) así como en su más conocida asociación en teoría de números a la distribución de los números primos. Estudiaremos con más detalle este último aspecto en la línea que comienza con los trabajos de Legendre y Tchebycheff, finalizando con Riemann.

\end{abstract}
\maketitle

\section{Introducci\'on}

La función logaritmo integral es conocida habitualmente por su conexión con la distribución de los números primos y con el "príncipe de las matemáticas", Carl F. Gauss (1777-1855). 

Sin embargo, objetivo importante del presente artículo es hacer visible que, además, el logaritmo integral tuvo un papel protagonista en el nacimiento de algunos conceptos y resultados relevantes en los albores del Análisis Matemático de finales del s. XVIII y principios del s. XIX.

Su puesta en escena vino de la mano de J. Bernoulli y L. Euler como parte de un problema de cálculo de primitivas. Sus particularidades llamaron la atención inicial de un buen número de matemáticos (Mascheroni, Plana, Bidone...). Del estudio detallado del logaritmo integral provienen los primeros cálculos precisos de la hoy conocida como constante de Euler-Mascheroni (ver Sección 2).

A finales del siglo XVIII la astronomía demandaba de las matemáticas nuevos métodos y herramientas para avanzar en su conocimiento, estando ambas disciplinas íntimamente relacionadas. Figuras tan relevantes como algunos de los protagonistas de nuestra historia estuvieron muy ligados a la astronomía: Olbers, Laplace, Legendre..., otros incluso fueron responsables de observatorios astronómicos como Bessel en Königsberg, Gauss en Götinga, Schumacher en Mannheim y Altona, Encke en Berlín...

Pues bien,  resolución de algunos problemas prácticos de astronomía planteados por Laplace, Kramp y otros, llevaban aparejadas ecuaciones que involucraban al logaritmo integral y que necesitaban del empleo de métodos numéricos de aproximación. En este contexto, en la sección 4, contamos el papel de J. Soldner que en 1809 bautiza esta función con el nombre y notación utilizado en la actualidad: el \textit{logaritmo integral} denotado como $li(x)$. Además desarrolló nuevos métodos de cálculo y construyó las primeras tablas de sus valores numéricos. 

En las secciones 5 y 6 emerge la figura del astrónomo y matemático Fiedrich W. Bessel (1784-1846).  Este era conocido entre otras cosas por ser el primero en calcular el paralaje de una estrella (61 Cygni) y por lo tanto su distancia a la Tierra. Bessel creó nuevos métodos que mejoraban los cálculos numéricos de la función $li(x)$, y lo hizo además motivado por complacer una necesidad de Gauss; si bien interés de éste era de una naturaleza muy diferente: la relación entre el logaritmo integral y la distribución de los números primos. 

En septiembre de 1810, \cite{Erm1852}, Bessel le cuenta a Olbers sus nuevos avances en el cálculo de $li(x)$ tal y como Gauss deseaba. Esta noticia ha sido reseñada por muchos autores, \cite{Lan1909}, \cite{Nar2000} y presentaba a Bessel como un simple observador-relator de los avances sobre los números primos de su colega Gauss. 

Sin embargo, al revisar las fuentes originales destacamos una carta \textit{anterior} de Bessel a Gauss\footnote{La fecha es de 26 de agosto de 1810. Fue publicada en 1880 \cite{Eng1880}, pero hasta donde sabemos no ha sido lo suficientemente reseñada.}, que muestra que la implicación de aquél fue de hecho de fundamental importancia para que Gauss pudiera reafirmar la estrecha relación entre la función $\pi(x)$, que cuenta los primos hasta un número $x$, y el logaritmo integral $li(x)$, es decir, lo que hoy se conoce como el Teorema del Número Primo (TNP). 

La colaboración que se establece entre ambos y que se realiza mediante intercambio de cartas dará lugar a que la función $li(x)$ aparezca en otros campos de la matemática aparentemente sin conexión. 
Así veremos en la sección 6 como el hoy conocido Teorema de Cauchy de variable compleja es anticipado por Gauss en 1811 y viene \textit{directamente} motivado por la pretensión de éste de hallar el valor de $li(a+bi)$. Años más tarde, Gauss introduce el famoso método de cuadratura que lleva su nombre y el \textit{único} ejemplo que proporciona es el cálculo de $li(200.000)-li(100.000)$. 

La conexión entre el logaritmo integral y la distribución de los números primos está bien documentada. A finales del siglo XVIII Legendre (Sección 3) y Gauss conjeturan de forma independiente el Teorema del Numero Primo, TNP. En la sección 7, relatamos detalladamente y cronológicamente las aportaciones que realizan Dirichlet, Tchebycheff, Hargreave... finalizando en los trabajos de B. Riemann, \cite{Rie1860} y \cite{Rie1800}. En éstos, la función $li(x)$\footnote{Riemann estimará \textit{numéricamente} tanto $li(x)$ como $\pi(x)$, como observaremos en sus papeles personales.} entrará en relación con la celebrada Hipótesis de Riemann. 
 
En la sección 8 aportamos una traducción comentada y contextualizada de la famosa carta de Gauss a Encke en la Navidad de 1849 en la que aquél da su particular visión de la historia de la conjetura de los números primos. En la última sección y como consecuencia de lo narrado en las secciones anteriores, hacemos incapié en las implicaciones que produce el desfase entre \textit{cuándo se publica} vs \textit{cuándo se obtiene} un resultado matemático. En ocasiones, una interpretación a posteriori no se ajusta a la realidad de los acontecimientos debido a que parte importante de los avances no fueron accesibles a todos los protagonistas implicados hasta mucho tiempo después. 

En el artículo hemos dado la máxima importancia a la búsqueda de información en las fuentes originales lo que se refleja en el alto número de citas bibliográficas utilizadas. Esto nos ha permitido mostrar las importantes aportaciones de algunos matemáticos que el paso del tiempo ha acabado difuminando. 

Como ayuda al lector concluimos la introducción aportando el siguiente cuadro resumen de los principales hitos que pretendemos exponer en este artículo.

\begin{center}
	\begin{tabular}{| c | c | c |}
		\hline
		\multicolumn{3}{ |c| }{Tabla resumen} \\ \hline
		Fecha del hallazgo & Hecho relevante & Año de publicación \\ \hline
		1737 & Identidad de Euler & 1744 \\	
		1763 & 1ª aparición de la integral, Euler & 1768 \\
		1790 & Mascheroni amplia investigaciones & 1790-2 \\
		1796, Mayo & Anotaciones de Gauss-conjetura TNP & 1917 \\
		1797 & Tabla de primos de Vega (400.000) & 1797 \\
		1798 & Fórmula de Legendre-conjetura TNP & 1798 \\
		1805 & Libro de Laplace	& 1805 \\
		1808 & Libro de Legendre, 2ª edic.  & 1808 \\
		1809 & Aparición de $li(x)$, Soldner & 1809 \\
		1810, 26 Agosto & Carta Bessel a Gauss & 1880 \\
		1810, 26 Agosto & Artículo de Bessel en Monastliche & 1810 \\
		1810, 1 Septiembre  & Carta Bessel a Olbers &	1852 \\
		1810-12	& Cartas entre Bessel-Gauss-Bessel & 1880 \\
		1811	& Tabla de primos de Chernac (1000000) & 1811 \\
		1811, 18 Diciembre &	Gauss anticipa Teorema Cauchy &	1880 \\
		1812 &	Tratado de Bessel &	1812 \\
		1814, 14 Septiembre	&	Método de cuadratura de Gauss & 1815 \\
		1838	&	Artículos de Dirichlet & 1838 \\
		1848, 24 Mayo &	Artículo de Tchebycheff	& 1848-52 \\
		1849, Julio & Artículo de Hargreave	&  1849 \\
		1849, 4 Diciembre & Carta de  Encke a Gauss	& 2018 \\
		1849, 24 Diciembre & Carta de Gauss a Encke	& 1863 \\
		1850, 7 Diciembre & Carta de Gauss a Dase &	1862 \\
		1859, Noviembre & Artículo de Riemann &	1860 \\ \hline
	\end{tabular}
\end{center}

\textbf{Notación:} La función $li(x)$ se define como: 
$\displaystyle li(x)=\int_{0}^x\frac{dt}{log\ t}, $ para valores $0\le x<1$ y  $\displaystyle li(x)=\lim_{\varepsilon\to 0^{+}}\left(\int_0^{1-\varepsilon}\frac{dt}{log\  t}+\int_0^{1+\varepsilon}\frac{dt}{log\ t}\right),$ si $x>1$. 

La función $\frac{1}{\log t}$ no es integrable en el entorno de $t=1$, motivo por el que se toma en la segunda expresión el llamado "valor principal". Esta definición recibe el nombre de "convención americana". La llamada "convención europea", ya presente en \cite{Lan1909}, toma como expresión $\displaystyle \int_{2}^{x}\frac{dt}{log\ t}=li(x)-li(2)$, ver una detallada justificación de ambas en \cite{Der2003}, Cap. 7. \ La razón de esta doble nomenclatura es de índole histórico y quedará reflejada a lo largo del artículo. En ocasiones, los matemáticos protagonistas de nuestra historia usarán la misma expresión $li(x)$ para representar objetos ligeramente distintos, lo que aclararemos en su momento.

\section{Antecedentes.}

Encontramos en los trabajos de Euler las primeras relaciones entre la todavía incipiente teor\'{\i}a de n\'umeros y el an\'alisis. En 1737 (\cite{Euler1744}, publicada en 1744)  Euler establece la identidad que lleva su nombre: 
$$\sum_{n=1}^{\infty}\frac{1}{n^s}=\prod_{p\ primo}\frac{1}{(1-\frac{1}{p^s})},\qquad s>1$$

A partir de ese momento, la identidad de Euler ser\'a para un buen n\'umero de autores la herramienta de arranque en la investigaci\'on de \textit{el problema la distribuci\'on de los n\'umeros primos en los n\'umeros naturales}.

\begin{figure}[htb]
	\includegraphics[scale=1.4]{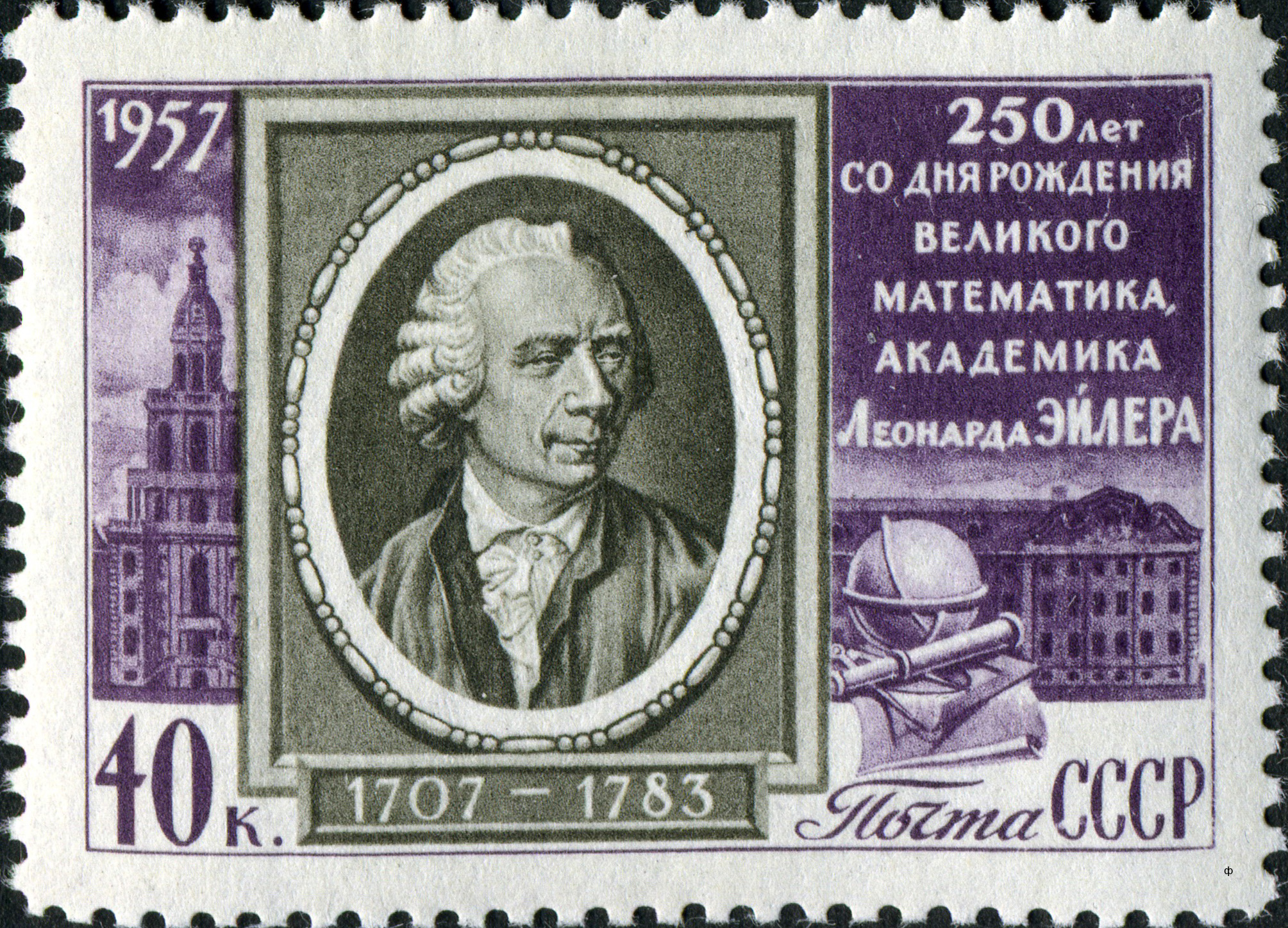}
	\caption{Leonhard Euler (1707-1783)}\quad
\end{figure}

Euler escribe y destaca su identidad tambi\'en en el caso $s=1$. A ojos de un lector actual, ese caso expresa la igualdad entre dos expresiones divergentes. Sin embargo, encierra mucha m\'as informaci\'on. Por un lado, le permite deducir la infinitud de los n\'umeros primos. Por otro, varios autores incluido el propio Euler se lanzaron a analizar de manera m\'as fina la suma/producto parcial de ambos lados de la igualdad. A modo de ejemplo,  cuando $N\to\infty$, Euler demuestra la convergencia de $ \Big(\sum_{n=1}^{N}\frac{1}{n}- \int_1^N\frac{dx}{x}\Big)\to \gamma\sim 0.57721...$ a la hoy llamada constante de Euler-Mascheroni, ver \cite{Euler1740} y \cite{Mas1790}.

Por otro lado tambi\'en se debe a Euler la primera aproximaci\'on a la expresi\'on protagonista\footnote{En esta secci\'on los autores citados usan la notaci\'on $\int\frac{dx}{\log x}$ para representar tanto una primitiva como una integral definida aunque sin especificar los extremos de integraci\'on, como har\'{\i}amos hoy en d\'{\i}a, puesto que \'estos se dan por entendidos por el contexto. Durante la secci\'on usaremos la notaci\'on de la \'epoca siempre que no de lugar a confusi\'on.} del presente art\'{\i}culo $\displaystyle \int\frac{dx}{\log x}$.

 Johann Bernoulli (1667-1748) recopiló y publicó sus escritos matem\'aticos en cuatro voluminosos tomos (Opera Omnia, 1742). En el tomo tercero considera el grupo de funciones $\int x^m (\log x)^\ell, m, \ell\in\N$ llegando por integraci\'on por partes, a la expresi\'on
$$\int x^m (\log x)^\ell dx=\frac{x^{m+1}}{m+1}  (\log  x)^\ell-\frac{\ell x^{m+1}}{(m+1)^2}\log  x^{\ell-1}+\frac{\ell(\ell-1)x^{m+1}}{(m+1)^3} (\log  x)^{\ell-2}-...$$

Citando a Bernoulli, Euler en 1768 \cite{Euler1768} toma como expresi\'on general a integrar $\displaystyle \int x^{m-1}(\log  x)^n dx, $ \textit{permitiendo a $n$ ser un entero negativo}. Mediante integraci\'on por partes y el cambio de variable $x^m=z$ el grupo de expresiones, para $n$ negativo, se reducen a $\displaystyle \int \frac{dx}{\log  x}$ \footnote{La notaci\'on del logaritmo cambia seg\'un los distintos autores. Euler por ejemplo utiliza $lx$. Para la comodidad de la lectura las hemos unificado en la notaci\'on moderna $\log x$.}. Euler dice de \'ella que \textit{hay que considerarla como una clase especial de funci\'on trascendental}, a la cual dedicar\'a una secci\'on un poco mas tarde.

Efectivamente, en el punto 228 mediante el cambio de variable $e^x=z$, Euler la transforma en la expresi\'on equivalente  $\displaystyle \int \frac{e^x dx}{x}$ y utiliza el desarrollo en serie de la exponencial para escribir las igualdades:
$$\int \frac{e^x dx}{x}=\gamma+\log x +\frac{x}{1}+ \frac{1}{2}.\frac{x^2}{1.2}+\frac{1}{3}.\frac{x^3}{1.2.3}+\frac{1}{4}.\frac{x^2}{1.2.3.4}+ etc $$
$$\int \frac{dz}{\log  z}=\gamma+\log\log z +\frac{\log z}{1}+ \frac{1}{2}.\frac{(\log z)^2}{1.2}+\frac{1}{3}.\frac{(\log z)^3}{1.2.3}+\frac{1}{4}.\frac{(\log z)^4}{1.2.3.4}+ etc  $$ 

Euler apunta a continuaci\'on los problemas de interpretaci\'on de las igualdades (sin entrar a resolverlos) que se presentan en los puntos $z=0$ y $z=1$ sentenciando que \textit{...la naturaleza de esta funci\'on trascendental se conoce poco}.

En 1790, L. Mascheroni \cite{Mas1790} sigue los pasos anunciados por Euler estudiando la expresi\'on $\int \frac{dz}{\log  z}$  para valores de $z\in (0,1)$. Como consecuencia de la precisi\'on del estudio, Marcheroni publica un valor aproximado de la constante $\gamma \simeq 0,577215 664901 532860 618112 090082 39$. Un total de 32 decimales de los cuales los primeros 19 decimales son correctos.

\begin{figure}[htb]
	\includegraphics[scale=0.5]{	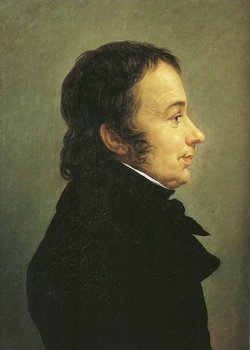}
	\caption{Lorenzo Mascheroni (1750-1800)}
\end{figure}

Finalmente, buena parte del contenido de las memorias aritméticas de Euler\footnote{Fueron recopiladas y publicadas por Fuss, en San Petersburgo, en 1849, en dos volúmenes.} se relaciona con la necesidad de ampliar las tablas de factores de n\'umeros, ver por ejemplo \cite{Euler1775}. El objetivo era estimular la ampliaci\'on de las tablas para poder progresar en las investigaciones. Es conocido a través de la correspondencia que Bernoulli animó a Lambert, Rosenthal y Felkel para que realizaran esta tarea. Euler no hará ninguna enumeración exhaustiva de números primos, pero sí será el primero en conjeturar que el n\'umero de primos hasta un valor $x$, denotado actualmente por $\pi(x)$, es aproximadamente igual a $x/\log x$, \cite{Nar2000}, \cite{Euler1764}. Como veremos más adelante, años más tarde el joven Gauss escribiría esta misma expresión en el margen de un texto de tablas numéricas, \cite{Schu1778}. 

La idea que parece subyacer en la mente de Euler y Gauss en estos años será expresada por G. H. Hardy en 1922 de la siguiente manera (y con lenguaje moderno): 

\textit{Hay infinitos números primos, su densidad decrece conforme avanzamos en la serie de los números naturales y tiende a 0 cuando tendemos a infinito. Más concretamente, el número de primos menores que $x$ es, en una primera aproximación $\frac{x}{log x}$. La probabilidad de que un número grande $n$ seleccionado al azar sea primo es aproximadamente alrededor de $\frac{1}{log n}$ y más precisamente, el logaritmo integral $li(x)=\int_2^x \frac{dt}{log t}$ da una muy buena aproximación al número de primos. \cite{Har1922}}

\section{Legendre, el hilo conductor}

En 1798 Legendre, \cite{Leg1798} (secciones 28-29,  en una nota a pie de página), escribe: \textit{Si hay $b$ números primos en la progresión $1,2,3...a$, para los siguientes valores de $a$ tenemos aproximadamente: $$a=10^1, \qquad	10^2 , 	\qquad	10^3, \qquad	10^4,\qquad	10^5\qquad …$$
$$\frac{b}{a}=1/2, \qquad 1/4, \qquad	1/6, \qquad	1/8 , \qquad	1/10 , \qquad …$$
de donde parece que podemos concluir $b=\frac{a}{2la}$}. Seguidamente propone 
 como aproximación de $\pi(x)$ la expresión $\frac{x}{A\log x+B}$, siendo $A$ y $B$ constantes y añade que \textit{"la determinación exacta de los coeficientes sería un problema curioso y digno de la atención de los analistas}.”

\begin{figure}[htb]
	\centering
	\includegraphics[scale=0.26]{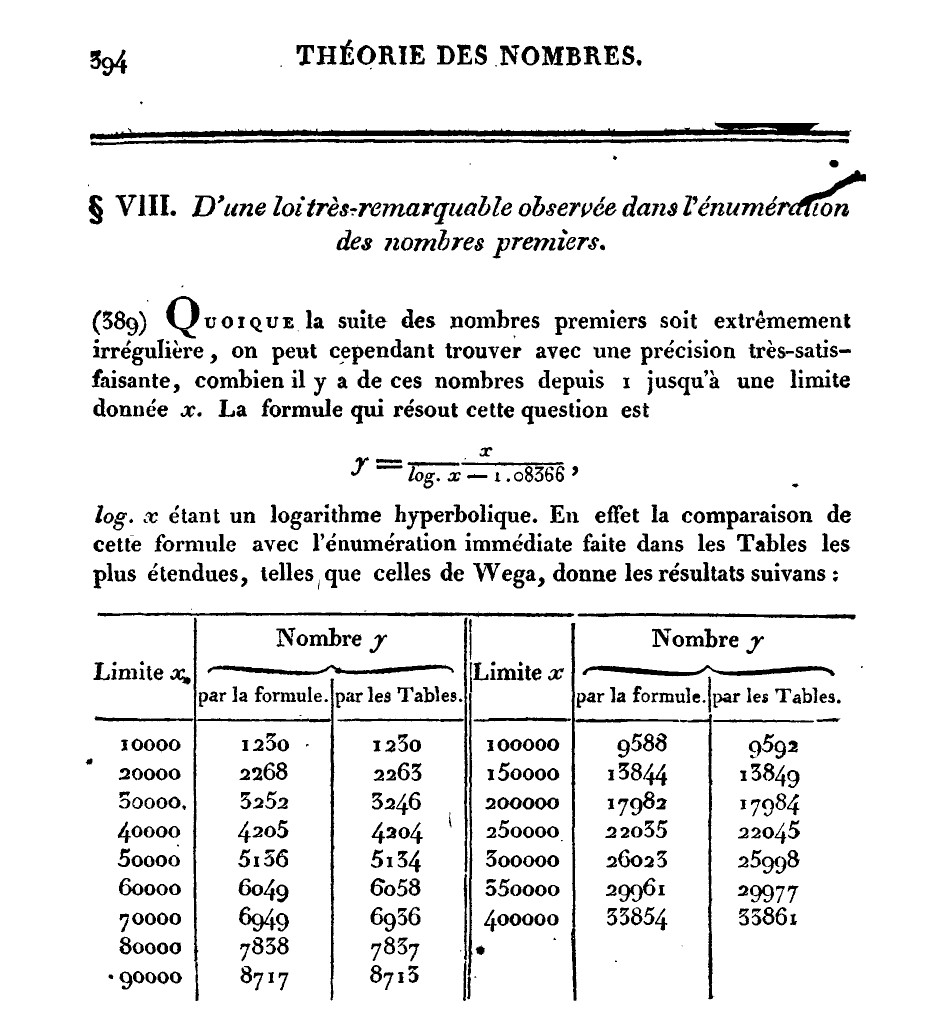}
	\caption{Legendre, 1808.}\label{fig:legendre}
\end{figure}

En la segunda edición  de 1808, \cite{Leg1808}, Legendre mejorará la expresión apoyándose en la identidad de Euler:

\textit{Aunque la secuencia de números primos es extremadamente irregular, podemos encontrar con precisión muy satisfactoria cuántos de estos números hay desde 1 hasta un límite dado $x$. La fórmula que resuelve esta cuestión es 
$$y=\frac{x}{\log x - 1.08366}$$
siendo $\log x$ un logaritmo hiperbólico."}

Mediante la inclusión de una tabla compara los resultados procedentes de su fórmula que, como Legendre comenta,  se ajustan bastante bien con las tablas de primos publicadas por Georgius Vega, \cite{Veg1797}. El número $1.08366$ llamará la atención del propio Gauss, que buscará sin éxito una interpretación del mismo, \cite{Gau1863}, \cite{Gol1973}.

Legendre no utilizará en este contexto la función $li(x)$. Sin embargo, sus artículos serán referenciados por la totalidad de los autores que, en el futuro, estudiarán la propiedad posiblemente más importante de $li(x)$, el de su relación con  problema de la distribución de los números primos.

\section{El Logaritmo Integral. J. Soldner y los primeros c\'alculos. }

En 1809 J. Soldner publica \textit{“Theorie et tables d’une nouvelle fonction transcendante”} \cite{Sol1809}. La funci\'on trascendente a la que se refiere dicho t\'{\i}tulo no es otra que la que denomina, por primera vez en la historia, con el nombre de \textit{Logaritmo Integral} y denota como $li(x)$. Nombre y notaci\'on ser\'an aceptados desde entonces. 

Soldner explica que su objetivo es el elaborar unas tablas num\'ericas de la $li(x)$ para facilitar la resoluci\'on de ciertos problemas f\'{\i}sicos tratados por P.S. Laplace \cite{Lap1805} y otros, que requer\'ian de gran precisi\'on en el c\'alculo. 

Asimismo, partiendo de las ideas de Euler y Mascheroni desarrollar\'a un nuevo m\'etodo que le permitir\'a calcular $li(x)$ para valores $x$ tanto pequeños como grandes: Soldner finalizó la tabla en $li(1280)$ al no tener necesidad de prolongarla. 

La expresi\'on principal en la que basa el c\'alculo para valores grandes de $x$ es: 

$$li(a+x)=li(a)+\frac{x}{\log a}-\frac{1.a A''}{1.2 (\log a)^2}y^2+\frac{2.a A'''}{1.2.3 (\log a)^3}y^3-\frac{3.a A''''}{1.2.3.4 (\log a)^4}y^4 + etc...$$
donde $y=\log(1+\frac{x}{a})$ y $A''=1,  A'''=1.A''-(\log a), A''''=2.A'''+(\log a)^2$ etc...

La f\'ormula permita calcular $li(a+x)$ a partir de $li(a)$, lo que obliga a que $\frac{x}{a}$ e $y$ sean pequeños. Con esto en mente, Soldner confecciona de manera recursiva su tabla (ver figura)

\begin{figure}[htb]
	\centering
	\begin{subfigure}{0.45\linewidth}
		\includegraphics[scale=0.6]{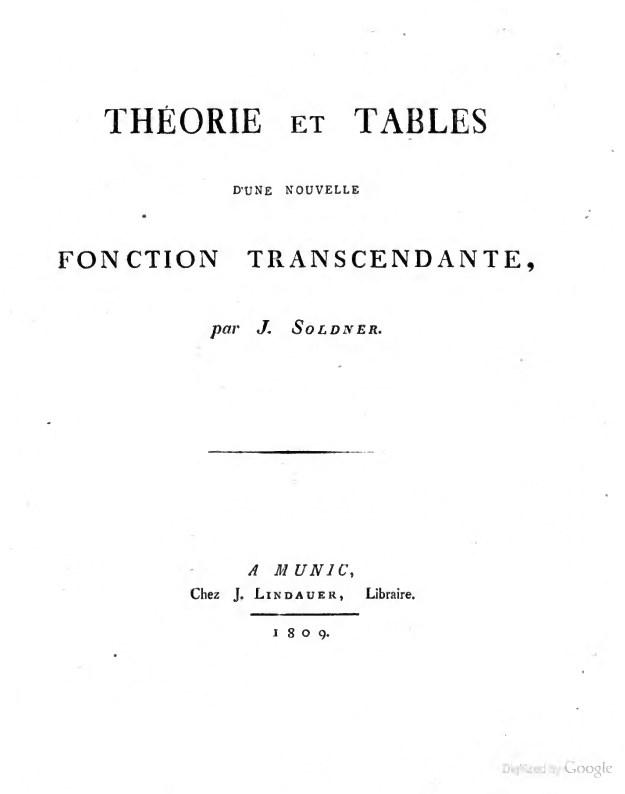}
	\end{subfigure}\quad\quad\quad
	\begin{subfigure}{0.45\linewidth}
		\includegraphics[scale=0.6]{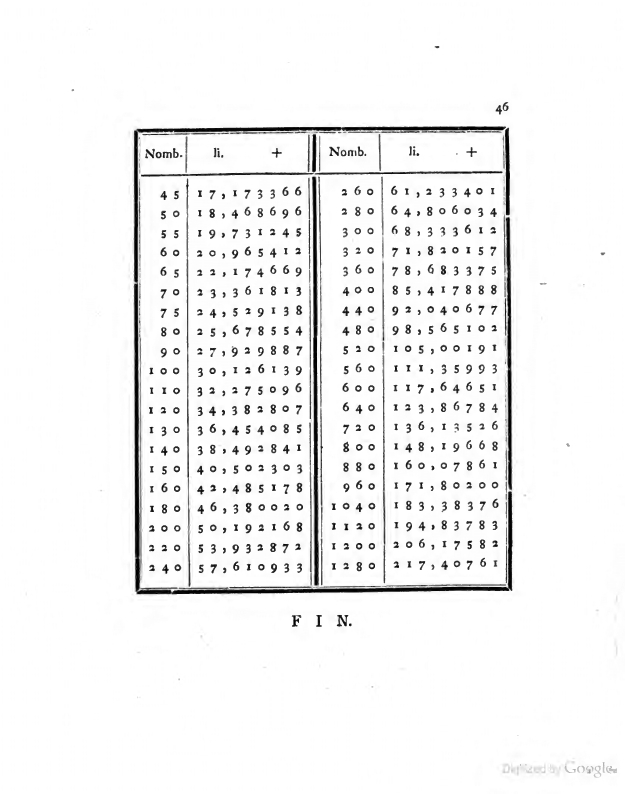}
	\end{subfigure}
	\caption{Tratado de Soldner sobre $li(x)$ de 1809}
\end{figure}

Soldner resuelve el problema de la no integrabilidad en el entorno de $t=1$ utilizando su valor principal. Igualmente, halla la \'unica ra\'{\i}z positiva de $li(x)$, la hoy llamada constante de Soldner-Ramanujan\footnote{Observamos que, por la propia definición de $\mu$, $li(x)= \int_{\mu}^x\frac{dt}{log\ t}$. Esta última expresión tiene la ventaja de evitar el problema de la no integrabilidad en $t=1$.}  y denotada por $\mu\sim 1,45136923488...$. Soldner de hecho escribe $li(1,4513692346)=0$ calculando correctamente sus 9 primeros decimales. Adem\'as calcula correctamente los primeros 22 decimales de la constante de Euler-Marcheroni, $\gamma\sim 0,5772156649015328606065..$.

El trabajo de Soldner no va a pasar desapercibido. En una carta fechada el 1 de septiembre de 1810 entre F. Bessel y H. Olbers publicada en 1852 \cite{Erm1852}, (ver la cita que de la misma realiza E. Landau \cite{Lan1909}; ver también \cite{Nar2000}) se puede leer: \bigskip

\textit{“Te doy noticias sobre una investigacion de la interesante integral $\int\frac{dn}{l n}$ y te digo algo que creo que te puede interesar. Soldner public\'o hace  un tiempo un tratado donde estudia este tema “Theorie et Tables d’une nouvelle Transcendante” del cual Gauss y Schumacher nos habían hablado y que Gauss había comentado que la integral hasta donde el sabe parecía estar conectada con el número de primos que hay hasta un n\'umero dado”}.  
\hfill Carta de Bessel a Olbers, 1-9-1810.\bigskip

Esta es la primera vez (documentada, aunque en comunicaci\'on privada) que el nombre de Gauss aparece vinculado a la funci\'on $li(x)$ y a los n\'umeros primos.

La referencia result\'o lo suficientemente sugerente a los autores como para buscar la carta original y comprobar de primera mano junto a sus protagonistas el contenido de la misma por si aportaba mas informaci\'on.  Para nuestra sorpresa, Bessel continuaba la carta con cuatro paginas de c\'alculos y explicaciones de c\'omo determinar $li(x)$ para valores \textit{muy grandes} de $x$.\bigskip

La figura de F. Bessel va a jugar un papel esencial, colaborando directamente con Gauss en el progreso de las investigaciones tanto sobre el logaritmo integral como en su relaci\'on con la distribuci\'on de los n\'umeros primos. 

\section{La aportaci\'on no contada de F. Bessel}

Friedrich Wilhelm Bessel (1784-1846) era considerado uno de los mejores calculistas del momento,  destacando por su minuciosidad en el c\'alculo de la \'orbita del cometa Halley as\'{\i} como por sus observaciones de las \'orbitas de Palas, Ceres y Juno. Por ello Olbers y Gauss (con el que comenz\'o su relaci\'on epistolar en 1804) le recomendaron como director del Observatorio de K\"onigsberg donde de hecho trabaj\'o a partir de 1810 y durante el resto de su vida. Adem\'as de como astr\'onomo, su nombre aparece ligado a multitud de aportaciones en matem\'aticas (funciones de Bessel, polinomios de Bessel...).

\begin{figure}[htb]
	\includegraphics[scale=0.8]{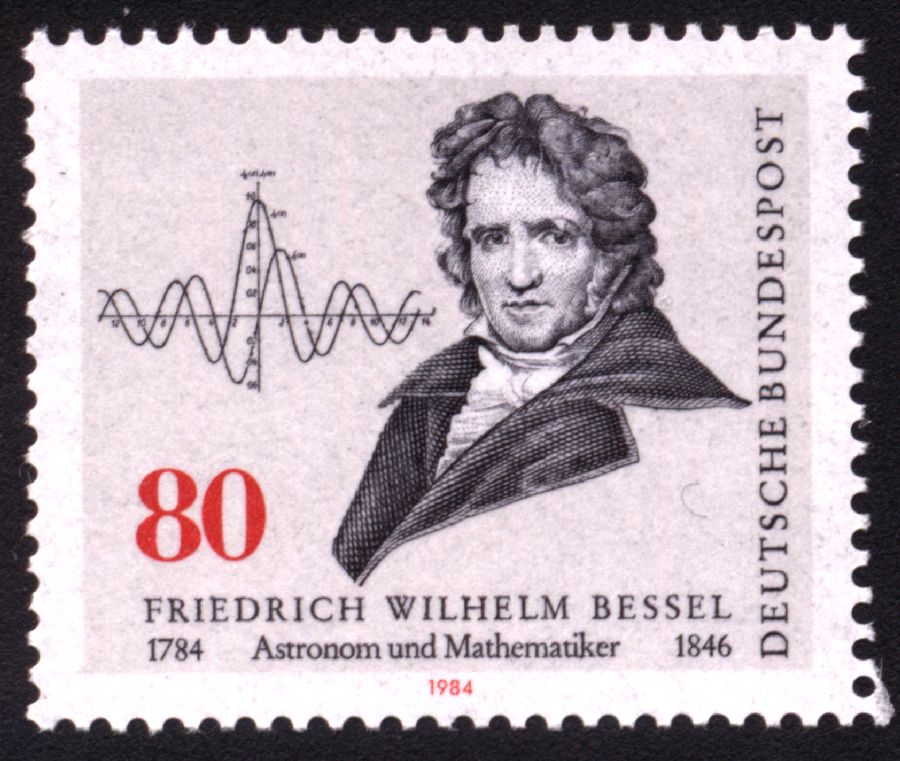}  
	\caption{Friedrich Wilhelm Bessel (1784-1846)}
\end{figure}

En la continuaci\'on de la mencionada carta Bessel sigue relatando a Olbers que ha conseguido desarrollar un nuevo m\'etodo que le permite...\bigskip

\textit{"...continuar con los calculos para $li(x)$ grandes \textbf{tal y como  Gauss desea} para poder confirmar su relaci\'on con los n\'umeros primos, ya que las tablas publicadas por Soldner solo llegan hasta x=1280. La integral tal y como la calcula Soldner no permite $x$ grandes sin acumular muchos errores..." }

\textit{"...le he dado a Gauss estas investigaciones con algo m\'as de detalle porque creo que aprecia este  tipo de desarrollos num\'ericos..." }

\hfill Carta de Bessel a Olbers, 1-9-1810.\bigskip

Esta nueva informaci\'on llev\'o a los autores a buscar alguna carta en fechas pr\'oximas a la carta del 1 de septiembre de 1810 por si Bessel y Gauss hubiesen intercambiado ideas, como se intuye por los comentarios realizados por Bessel a Olbers. La sorpresa lleg\'o con una carta (que hasta donde sabemos no ha sido reseñada en la historiograf\'{\i}a) fechada el 26 de agosto y que ademas coincid\'{\i}a en el tiempo (ver Tabla Resumen) con un art\'{\i}culo enviado por Bessel a Lindenau y publicado en Monatliche Correspondenz con el t\'{\i}tulo \textit{Ueber das Integral $\int\frac{dx}{l x}$}, \cite{Bes1810}. Lo que claramente indicaba que las investigaciones hab\'{\i}an comenzado con anterioridad.

A continuaci\'on presentamos el papel de Bessel en estos primeros pasos traduciendo un extracto de la carta del 26 de Agosto de 1810: \medskip

\textit{"...comparto con usted el resultado de estos hallazgos ya que puedo suponer que le interesan. La investigaci\'on de la que hablo es sobre la integral $\int\frac{dx}{l x}$ y fue escrito por Soldner y que recib\'{\i} hace unas semanas. Sin duda conoces este pequeño trabajo y las dificultades encontradas por el autor que no le permitieron extender la tabla para la funci\'on $\int\frac{dx}{l x}$ que el denomina $li\ x$ (logaritmo-integral de x) para valores mas alla de $x=1280$. El camino que ha tomado hasta ahora no difiere mucho del proceso de aplicar el teorema de Taylor y es bastante tedioso, de hecho inadmisible, ya que no da ningun control sobre los c\'alculos y no puede revelar los errores que podr\'{\i}an haberse acumulado.  Por lo tanto, busqu\'e otro m\'etodo para hacer la integral y encontr\'e lo siguiente".}

Bessel utiliza un ingenioso artificio, que le ocupa unas seis p\'aginas, para representar $li (\frac{x}{a})$ (para $\frac{x}{a}>1$) mediante:
$$li (\frac{x}{a})=li(x)+x\Big\{\frac{A'}{lx}+\frac{A''}{(lx)^2}+\frac{A'''}{(lx)^3}+etc...\Big\}$$
donde $A'=\frac{1}{a}-1, A''=A'+\frac{\log(a)}{a}, A'''=2A''+\frac{(\log(a))^2}{a}, A''''=3A'''+\frac{(\log(a))^3}{a}...$. A continuaci\'on se interesa en los valores $x=10^n$ y $a=10$: 

\textit{"se adivina facilmente, mi querido amigo, porque te estoy dando estos desarrollos una vez expresado el deseo de conocer la $li(x)$ para valores muy grandes de la $x$, para poder alabar el hermoso comentario de la conexión con los primos. Entonces tengo el valor de $li(x)$  de 1 a 1.000.000 que he calculado y as\'{\i} podr\'{\i}a continuar f\'acilmente ese c\'alculo en gran medida si tuviera alg\'un inter\'es. Para avanzar r\'apido he tomado a=10 y $x=10^n$..."}

La expresi\'on simplificada obtenida le permitir\'a ir calculando $li(x)$ para potencias de 10

$$li ( 10^n)= li (10^{n-1})- \frac{A'}{\log 10}\frac{10^n}{n}-\frac{A''}{(\log 10)^2}\frac{10^n}{n^2}-\frac{A'''}{(\log 10)^3}\frac{10^n}{n^3}+...$$

y comenta satisfecho las buenas estimaciones de los errores de aproximaci\'on que se cometen:

\textit{"he calculado tantos coeficientes que siempre puedo usar la serie, incluso si no se supone que $\frac{x}{10}$ sea mucho mayor que 1; por lo general, si x es un poco grande, de 8 a 10 terminos son suficientes, por ejemplo con $10^5$, el termino 11 agrega 10 decimales correctos"}

Finalmente, proporciona a Gauss una serie de tablas con los coeficientes calculados y la tabla de los valores de la $li(x)$ comparada con el de n\'umero de primos hasta $x$. \bigskip

\begin{center}
	\begin{tabular}{| c | c | c | c | }
		\hline
		\multicolumn{4}{ |c| }{Tabla de Bessel, 1810} \\ \hline
		Funci\'on $li$ & Valor & N\'umero de primos & Exceso \\ \hline  
		li 1000 & 177,609655  & 169 & + 8,61  \\
		li 10000 & 1246,137247 & 1230 & + 16,14  \\
		li 100000 & 9629,809041 & 9593 & + 36,81  \\
		li 200000 & 18036,052159 & 17983 & + 53,05 \\
		li 300000 & 26080,215589 & 25997 & + 83,21 \\
		li 400000 & 33922,621995 & 33859 & + 63,62 \\
		li 1000000 & 78627,549277 & --- & --- \\ \hline
	\end{tabular}
	\label{Bessel}
\end{center}
\bigskip

Hay que resaltar que la columna de n\'umeros primos, como la publicada por Legendre, tan s\'olo llega hasta 400.000. En ambos casos, los datos hab\'{\i}an sido tomados de las tablas de primos de G. Vega \cite{Veg1797}.


Bessel corrige algunos datos de la tabla de Vega, incluye en el recuento el 1 y el 2 y se lamenta de que las unicas tablas existentes y confiables sean las de Felkel y Vega y 

\textit{"que si supieran el curso posterior de los primos, seria f\'acil encontrar una funci\'on que se adapte mejor y que los conectara"}

Bessel sugiere varias alternativas a la $li(x)$ como \textit{"las funciones  $\sum\frac{1}{lx}$ , o de forma mas general
$\sum\frac{1}{l(x+a)}$, con a constante; o con la integral $\int\frac{dx}{l(x+\frac{1}{x})}$ ,   o tal vez $\sum\frac{1}{l(x+\frac{1}{x})}$, etc". }


Hay que hacer notar que las primeras referencias sobre la conexi\'on de la $li(x)$ con los números primos hasta un valor $x$, tan s\'olo se encuentran en las comunicaciones  epistolares de Bessel con Olbers y Gauss. Las publicaciones de Bessel sobre el tema \cite{Bes1810}, \cite{Bes1811}, \cite{Bes1812} (ver las obras completas de Bessel en \cite{Eng1876}), nunca hacen referencia a dicha conexi\'on. 
En nuestra opini\'on es esta la principal raz\'on por la que el nombre de Bessel no aparece en la historia sobre la distribuci\'on de los primos. Como hemos visto, proporciona a Gauss las primeras estimaciones fiables y conocidas de $li(x)$ hasta 1.000.000, y de hecho le proporciona un m\'etodo que permite calcularla para valores mucho mayores de $x$. 

A diferencia del m\'etodo de Soldner que relaciona $li(x+a)$ con $li(x)$, el de Bessel conecta $li(ax)$ con $li(x)$. Un salto cualitativo (de un paso aritm\'etico a otro geom\'etrico) fundamental. Otros estudios de la función $li(x)$ en la primera mitad del siglo XIX aparecen en Bidone \cite{Bid1805}, Spence \cite{Spe1809}, Lacroix \cite{Lac1819}, Plana \cite{Pla1821}, Bretschneider \cite{Bre1837}, de Morgan \cite{Mor1842}, Schl\"omomilch \cite{Sch1847}.


La respuesta de Gauss llegará el 21 de octubre de 1810, \cite{Eng1880}: 

\textit{"Me has dado una gran alegr\'{\i}a con tu investigaci\'on de la integral $\int\frac{dx}{log(x)}$. B\'asicamente hab\'{\i}a tomado el mismo camino, pero no hab\'{\i}a llegado tan lejos".}

 Gauss prefiere expresar la integral a trav\'es de un cambio de variable como $\int\frac{e^y\, dy}{y}$ y comunica a Bessel sus resultados iniciales.
 
 La relaci\'on epistolar entre ambos sobre la $li(x)$ continuar\'a durante dos años. Bessel acabar\'a publicando en 1812 un extenso tratado sobre el logaritmo-integral, \cite{Bes1812}, a cuyo borrador Gauss tendr\'a acceso. En él se encuentran detalladas las aportaciones que había comunicado a Gauss por carta, ampliadas con nuevas tablas numéricas tanto de $li(x)$ como de sus incrementos $li(x)-li(x')$. Además añade nuevos desarrollos en serie, en fracción continua, así como relaciones con otras funciones trascendentes. En el artículo 18, \cite{Bes1812}, Bessel escribe una relación que provocar\'a en Gauss una sorprendente e interesante respuesta. Las expresiones son 
 $$\int\frac{\sin x}{x}dx=\frac{li(e^{ix})-li(e^{-ix})}{2ix}\qquad \qquad 
 \int\frac{\cos x}{x}dx=\frac{li(e^{ix})+li(e^{-ix})}{2x}$$

\section{Los descubrimientos no publicados de Gauss}

En la sección anterior hemos relatado la primera vez (documentada) en la que el nombre de Carl Friedrich Gauss (1777-1855) aparece relacionado con la función $li(x)$ y la distribución de los números primos. Aunque no public\'o nada al respecto, \cite{Lan1909}, \cite{Nar2000}, Gauss, como queda constancia por la recopilaci\'on p\'ostuma de sus obras completas, \cite{Gauss}, se sinti\'o atraido por los números primos durante toda su vida.

\begin{figure}[htb]
	\includegraphics[scale=0.08]{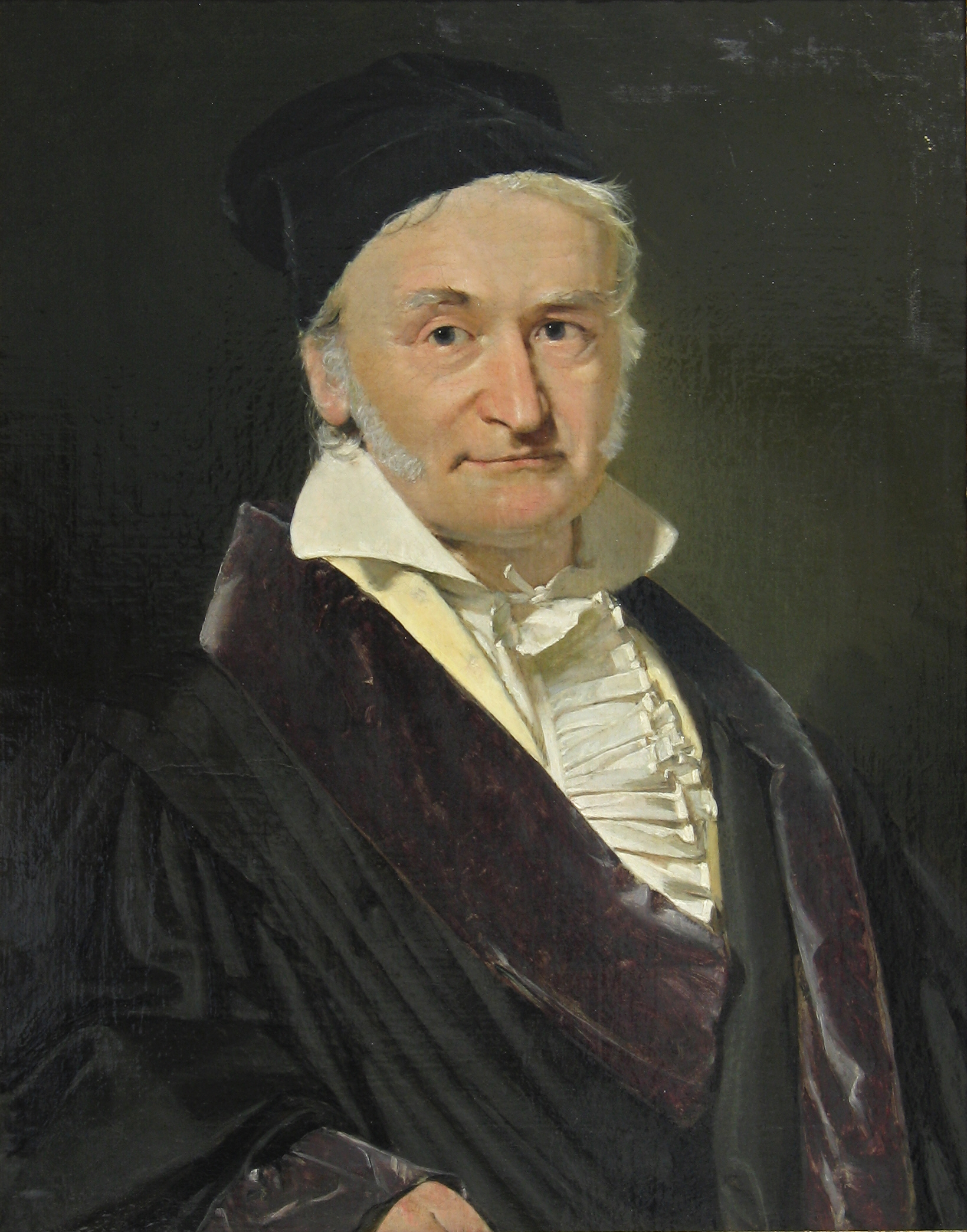}  \quad \quad \quad \quad
	\caption{Johann Carl Friedrich Gauss (1777-1855)}\quad
\end{figure}

Siguiendo el hilo temporal, en una segunda carta fechada el 19 de octubre de 1810, \cite{Eng1880}, Bessel comunica a Gauss nuevos avances que permiten c\'alculos m\'as directos de $li(x)$. La respuesta de Gauss a dicha carta de Bessel así como al borrador recibido se produce el 18 de diciembre de 1811, \cite{Eng1880}. 

 Gauss comienza manifestando su felicidad por haber recibido finalmente el borrador del tratado de Bessel (\cite{Bes1812})  y que ahora tiene por fin los valores de la $li(x)$ para grandes cantidades, adem\'as se congratula de haber recibido la tabla de n\'umeros primos de Chernac hasta 1,020,000 recientemente publicada.

Aqu\'{\i} Gauss revela:

\textit{"....he venido calculando progresivamente primos de miriada en miriada\footnote{"De miriada en miriada" significa "de 10.000 en 10.000"} y as\'{\i} compararlos con el valor de la integral $\int\frac{e^x-1}{ x}dx$".} 

\textit{"Me has hecho darme cuenta de que me gustar\'{\i}a mostrar estos resultados en nuestros anuncios acad\'emicos: nuestra amistad, querid\'{\i}simo Bessel, hace que sea mi deber, antes de hacerlo, tener una conversaci\'on por escrito contigo sobre un punto u otro, donde mi punto de vista  no esta del todo de acuerdo con el tuyo. Entonces toma los siguientes comentarios a bien y comp\'artelos conmigo con toda la franqueza y tan abiertamente como yo... compartimos un principio com\'un para los dos "no hay verdaderas controversias en matem\'aticas", por lo que nos entenderemos a trav\'es del intercambio mutuo de nuestras ideas..."}

Aqu\'{\i} Gauss da el gran salto y plantea a Bessel que la funci\'on no s\'olo se aplique para los numeros reales grandes, d\'andose cuenta de que\textit{ tambi\'en se puede aplicar para los n\'umeros imaginarios con los mismos derechos y sin restricciones}.\bigskip

 \textit{"Creo que debo asumir que esencialmente est\'a de acuerdo conmigo en este punto y que sus explicaciones en el art. 18, ya muestran que de ninguna manera desea bloquear su camino a contemplar la función $li(a+bi)$"}\bigskip

Es en este momento en que Gauss introduce por primera vez el problema de la integraci\'on de funciones complejas pregunt\'andose: \bigskip

\textit{"?`Qu\'e debemos entender por $\int \Phi(x) dx$ cuando $x=a+bi$?... Afirmo ahora que la integral $\int \Phi(x) dx$ tiene un valor \'unico a\'un tomada sobra varias trayectorias siempre que $\Phi(x)$ tome un valor \'unico y no se haga infinita en el espacio comprendido por las dos trayectorias. Este es un teorema muy bello cuya demostraci\'on no es dif\'{\i}cil y que proporcionar\'e en otro momento"}. \bigskip

Es decir, Gauss se adelanta en más de una década a lo que, desde 1825, ser\'a conocido como el Teorema de Cauchy.  Varios autores, \cite{Kli1972}, han destacado este hallazgo tan importante en el campo de la variable compleja. Lo que ha pasado más inadvertido es que una importante fuente motivadora del descubrimiento es justamente la funci\'on $li(x)$ para valores complejos. 

Seguidamente Gauss alerta sobre el efecto que produce en la integral de $1/x$ el que la funci\'on logaritmo sea multivaluada, dependiendo del camino de integraci\'on. Nos est\'a introduciendo en el concepto de índice. En sus propias palabras: \bigskip

 \textit{"...si define usted $log(x)$ a traves de $\int \frac{1}{x}$, desde $x=1$, puedes llegar a $log(x)$ sin incluir el punto $x=0$; o si pasas por el mismo una o más veces, cada vez se agrega la constante $+2\pi i$ ó $-2\pi i$: por lo tanto, los múltiples logaritmos de cada número son muy claros". }\bigskip

Gauss piensa que \textit{"extender las investigaciones a argumentos imaginarios dará lugar a resultados muy interesantes"} y se decanta finalmente (como ha dicho al inicio de la carta) por  $\int\frac{e^x -1}{x}dx$ \textit{"porque sospecho que la primera dará resultados mejores"}...\textit{"Yo preferiria, en vez de elegir la funci\'on del Sr. Soldner  $li( x)=\int\frac{dx}{log x}$, introducir una nueva función porque una función de una sola forma siempre puede considerarse más clásica y simple que una multivaluada, especialmente porque $log(x)$ es una funcion multivaluada. Podría ser bueno presentar $\int \frac{e^x-1}{x}dx$, o mejor $\int\frac{e^x}{x}dx$... que llamaré brevemente $Eix$, o exponencial-logaritmica, con su propio carácter y nombre especialmente porque las tareas de la física que conducen a $li(x)$ son generalmente exponenciales de $x$ en sí mismas."}\bigskip

La expresi\'on sugerida por Gauss,` $Ei(x)=\lim_{\varepsilon\to 0}\Big(\int_{-\infty}^{-\varepsilon}\frac{e^t}{t}dt+\int_{\varepsilon}^{x}\frac{e^t}{t}dt\Big)$ recibe en la actualidad el nombre de "integral exponencial" y se tiene $Ei(x)=li(e^x)$,  \cite{Gla1882}.

Gauss finaliza la carta refiriéndose a su propio c\'alculo de la constante de Euler-Mascheroni: 
\bigskip

\textit{"Hace mucho tiempo calculé la constante ... con 23 cifras...y la encontré diferente de la de Mascheroni a partir del decimal 20 ... creo que mi número es el correcto, pero yo no guardo las cuentas."} \smallskip

A lo largo de 1812 continuar\'a el intercambio de ideas por carta sobre la comprensi\'on de las expresiones en cuanto a caminos de integraci\'on, convergencia de desarrollos en serie, valores imaginarios etc... A modo ilustrativo, el 5 mayo 1812 Gauss comenta a Bessel las últimas novedades de lo que será su trabajo sobre las funciones hipergeométricas (de las que $li(x)$ es un caso particular), \cite{Gau1813}.

Sin embargo, a pesar de la insistencia de Gauss, Bessel no profundizará en las cuestiones de variable compleja que tanto interesa a aquel \textit{"lamento estar tan lejos de ti y espero conocer tu tratado pronto..."} a la vez que espera \textit{"que pronto me diga el resultado de su hallazgo de los n\'umeros primos que estoy muy intesesado en ver crecer un nuevo triunfo para usted"}.\footnote{Registramos tres cartas de Bessel a Gauss referentes al tema: 12 de enero, 26 de marzo y 27 de noviembre y dos de Gauss a Bessel: 5 de mayo y 31 de agosto, \cite{Eng1880}}\bigskip

El inter\'es de Gauss en la $li(x)$ parece no desvanecerse. En 1815 publica su \textit{“Methodus nova integralium valores per approximationem inveniendi”}, \cite{Gau1815}, en el que introduce las reglas de cuadratura que llevan su nombre. Estas le van a permitir aproximar valores de integrales (definidas) de una forma mucho más precisa y eficaz con notable ahorro de cálculo. Al final del tratado ilustra su nuevo m\'etodo con un \'unico ejemplo.  El ejemplo elegido, curiosamente, es el cálculo de la integral de la forma $li(x+h)-li(x)=\int_{x}^{x+h}\frac{dx}{log x}$, concretamente en el intervalo que va de 100,000 a 200,000. Gauss lo compara con el obtenido por Bessel, \cite{Eng1876}, \cite{Bes1812}, verificando así su eficacia. Esta ser\'a la única ocasión \textit{publicada en vida} de Gauss donde su nombre aparece junto a la funci\'on $\int\frac{dx}{log x}$ y al nombre de Bessel (sin aludir en ningún momento a su relación con los números primos). 

\section{La relaci\'on de $li(x)$ con los n\'umeros primos}

En esta sección explicamos el más conocido uso de la función $li(x)$: el de su relación con la distribución de los números primos. 
En un orden cronológico daremos un repaso a los artículos publicados sobre esta cuestión, para terminar con el célebre artículo de B. Riemann de 1859. 

Recordemos que en las dos secciones anteriores las discusiones alrededor de dicha relación tan solo fueron epistolares, sin trascender a publicaciones científicas. 

Por tanto, desde el punto de vista \textit{"oficial"}, esta sección comienza después de los trabajos de Legendre cuyo "Theorie des Nombres" es reeditado en 1808 y 1830. Esta última edición se convertirá en el punto de partida para todos los autores que narramos a continuación. 

En varios artículos fechados alrededor de 1838, \cite{Dir1889}, G.L. Dirichlet (1805-1859) plantea la búsqueda de leyes asintóticas para varias funciones en teoría de números. Sus trabajos demuestran, en particular, la existencia de infinitos números primos en progresiones aritméticas. Afirma que con sus métodos puede demostrar la fórmula propuesta por Legendre $\pi(x)\sim \frac{x}{\log x}$ y que la publicará proximamente, algo que nunca hizo. Hacemos notar que uno de sus artículos de 1838, que envía a Gauss, contiene una nota manuscrita en la que propone como primera aproximación de $\pi(x)$ la expresión $\sum_{n\le x}\frac{1}{log n}$, es decir la versión discreta de $li(x)$. 

El siguiente momento relevante aparece con los trabajos de P.L. Tchebycheff (1821-1894). Desde su llegada en 1847 a S. Petersburgo desde Moscú, éste se involucra en la recopilación y publicación de los trabajos de Euler. Motivado por ello, comienza investigar en teoría de números. Su principal trabajo es leído el 24 de mayo de 1848 ante la Academia de St. Petersburgo formando parte de su tesis doctoral "Teoría de congruencias" \cite{Tch1848}, \cite{Tch1899}. Tchebycheff muestra conocer las publicaciones de Legendre y Dirichlet. 

Tchebycheff considera el resultado empírico de Legendre observando que se ajusta muy bien al compararlo con las tablas de primos publicadas en su “Theorie des nombres” de 1830. Sin embargo su investigación le lleva a apuntar hacia la función  $\int_2^x\frac{dt}{log t}$ como la más apropiada\footnote{El uso de esta expresión $li(x)-li(2)$ llevará a la llamada "convención europea", utilizada entre otros por Landau, \cite{Lan1909}.}. Estudiará en detalle la diferencia $\pi(x)-\int_2^x\frac{dt}{log t}$, concluyendo en particular que la expresión de Legendre es incorrecta.
Asimismo es consciente de la limitación de las tablas disponibles y predice que se tendrá que ir a valores de $x$ superiores a 10 millones para que sea apreciable su mejora, lo que no era entonces comprobable.

\textit{"Es fácil convencerse por las tablas de números primos que la integral $\int_2^x\frac{dt}{log t}$ para $x$ grandes, expresa con suficiente precisión cuantos números primos hay inferiores a $x$. Estas tablas son demasiado pequeñas para poder decidir sobre la superioridad de la fórmula $\int_2^x\frac{dt}{log t}$ sobre la de Legendre $\frac{x}{log x - 1,08366}$".  }

\begin{figure}[htb]
	\centering
	\begin{subfigure}{0.45\linewidth}
	\includegraphics[scale=0.13]{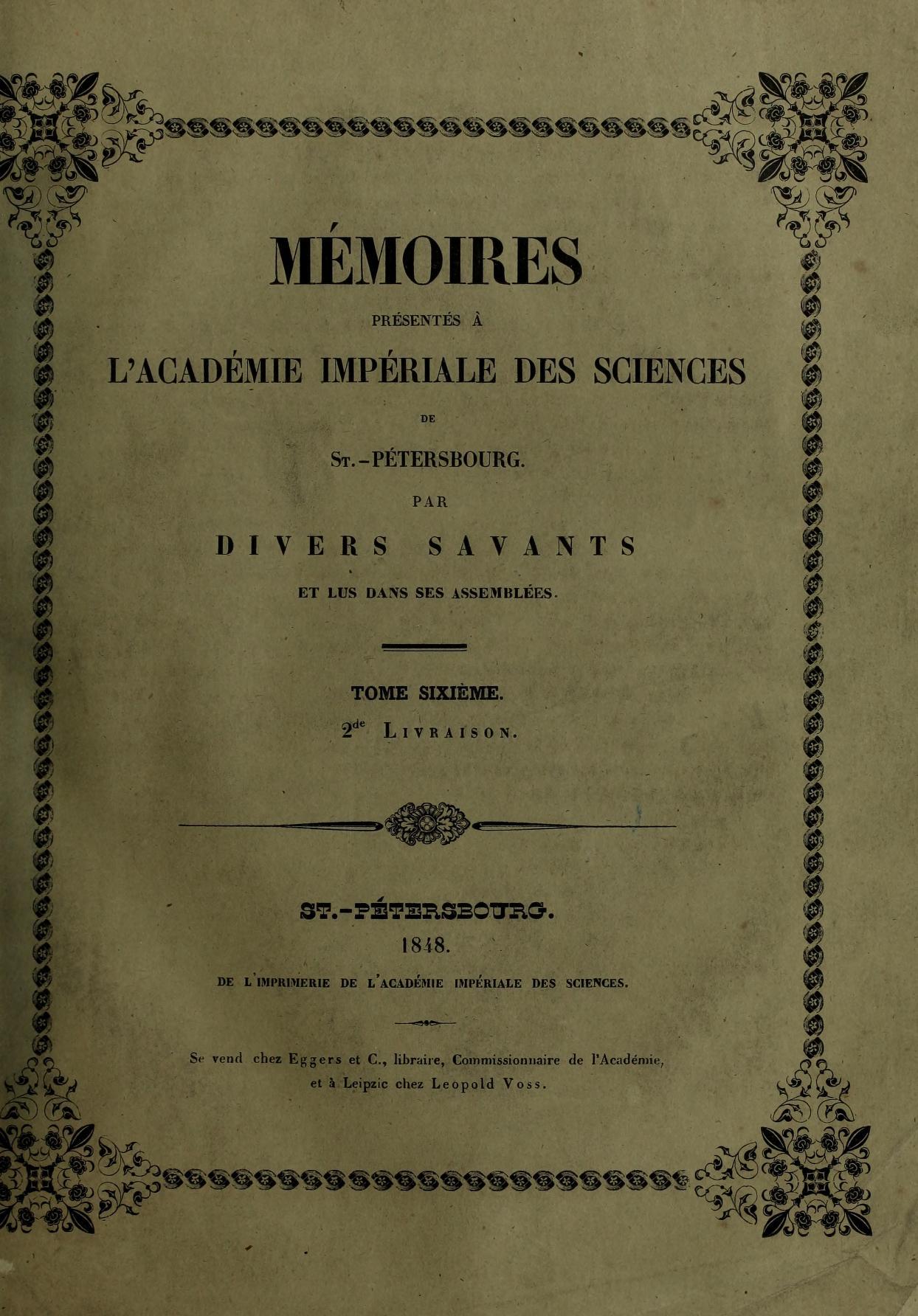}
	\end{subfigure}\quad\quad\quad
	\begin{subfigure}{0.45\linewidth}
		\includegraphics[scale=0.13]{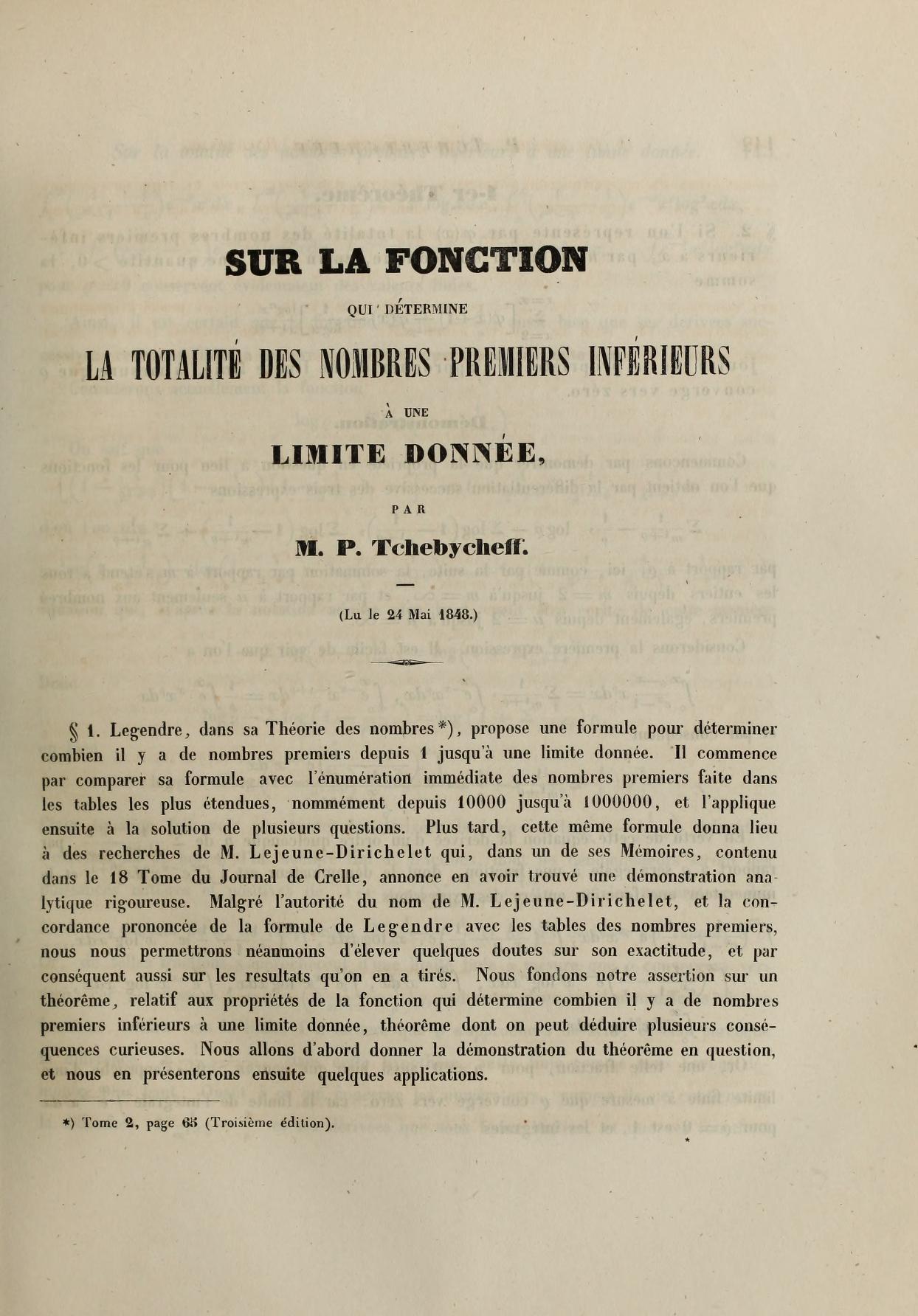}
	\end{subfigure}
\caption{Artículo de Tchebycheff de 1848}
\end{figure}

\medskip
La importancia del trabajo de Tchebycheff será enorme. Se convertirá en el punto de partida de todas las investigaciones posteriores sobre el problema de la distribución de los números primos.

\bigskip

La revista "Philosophical Magazine" de julio de 1849, contiene un art\'{\i}culo del abogado Sr. C.J. Hargreave, titulado \textit{Analytical Researches concerning Numbers},\cite{Har1849}. Hargreave es conocedor de los trabajos de Legendre y de Soldner a través del texto de de Morgan \cite{Mor1842}, cuyo texto es la principal referencia de la época en lengua inglesa.    

Aunque sus aportaciones no fueron de trascendencia, propone argumentos heurísticos más que verdaderas demostraciones, este artículo de Hargreave es considerado como el primer momento en que aparece \textit{publicada} la conexión entre la distribución de los números primos y la función $li(x)$, \cite{Nar2000}. 

\textit{"En este artículo tendré la ocasión de emplear la conocida integral trascendental 
$\int_0^x\frac{dt}{log t}$... bautizada con el nombre de logaritmo-integral denotada por $li(x)$ y que ha sido tabulada por Soldner..."}

En un segundo artículo, \cite{Har1854}, estimó el número de primos ($\pi(x)$) inferior a un millón, a cinco millones y a diez millones utilizando el método de conteo de Legendre,  \cite{Leg1808}, comparándola con $li(x)$.

En los años siguientes a 1850, Tchebycheff presentaría a la Academia otros trabajos, \cite{Tch1899}, en los que profundizaría en el problema de la distribución de los números primos demostrando el postulado de Bertrand sobre la existencia de (al menos) un primo entre $n$ y $2n$ y estimando de manera más precisa la función $\pi(x)$ $$0.92129...\le\frac{\pi(x)}{x/\log x}\le 1.105 55...$$

Asimismo, probó que si el límite de $\pi(x)/li(x)$ cuando $x\to\infty$ existe, éste debe ser igual a 1. Los artículos fueron reimpresos en la entonces conocida como Liouville Journal en 1852. 

\subsection*{La importancia de las tablas de factores}

En 1850, Gauss envía una carta, \cite{Das1862}, al calculista y matemático alemán Z. Dase (1824-1861), que era conocido por su gran talento para la realización de tablas matemáticas, instándole a ampliar las tablas existentes de números primos. La carta comienza justificando la utilidad de las mismas:

\textit{"La necesidad de descomponer los números existentes en sus factores, o de reconocer la imposibilidad de separarlos, se presenta muy a menudo a todos los que tienen que lidiar mucho con cálculos numéricos. Pero este esfuerzo aumenta a medida que aumentan los números, de modo que incluso un calculista muy experimentado puede necesitar horas, incluso días completos, para un solo número... Pero si las tablas se han calculado de una vez por todas, dan la respuesta en la medida de lo posible, el número puede ser pequeño o grande..."
}\bigskip

\hfill Carta de Gauss a Dase, 7-12-1850.\bigskip

Las tablas de números primos mas conocidas hasta finales del siglo XVIII eran las publicadas por J.H. Lambert (1770,1798),\cite{Lam 1770}, J.C. Schulze (1778), \cite{Schu1778}, y G. Vega (1797), \cite{Veg1797}, siendo esta última la más completa alcanzaba solamente hasta el primo 400031. Gracias a los archivos personales, hoy sabemos que Gauss poseía todas éllas, siendo alguna obsequio de sus mentores, \cite{Gau2005}. En la segunda decada del siglo XIX, las tablas de L. Chernac (1811), \cite{Che1811}, alcanzarían el listado de primos hasta el millón y las de J.C. Burckhardt (1814-17), \cite{Bur1814}, llegarían hasta los 3 millones. El propio Gauss junto con su ayudante Goldschmidt elaboraron tablas propias hasta los 3 millones como queda constancia en sus Werke \cite{Gauss}. 

Por otro lado, A. Crelle entre 1830-40 trabajó en la elaboracion de tablas y completó el 4º, 5º y 6º millon depositando en 1842 en la Academia de Berlin un manuscrito con los numeros primos hasta 6 millones que nunca fue publicado en contra del deseo de Gauss expresado en carta a Schumacher del 25 de enero de 1842, \cite{Sch1860}.

 \begin{figure}[htb]
	\centering
	\begin{subfigure}{0.45\linewidth}
		\includegraphics[scale=0.046]{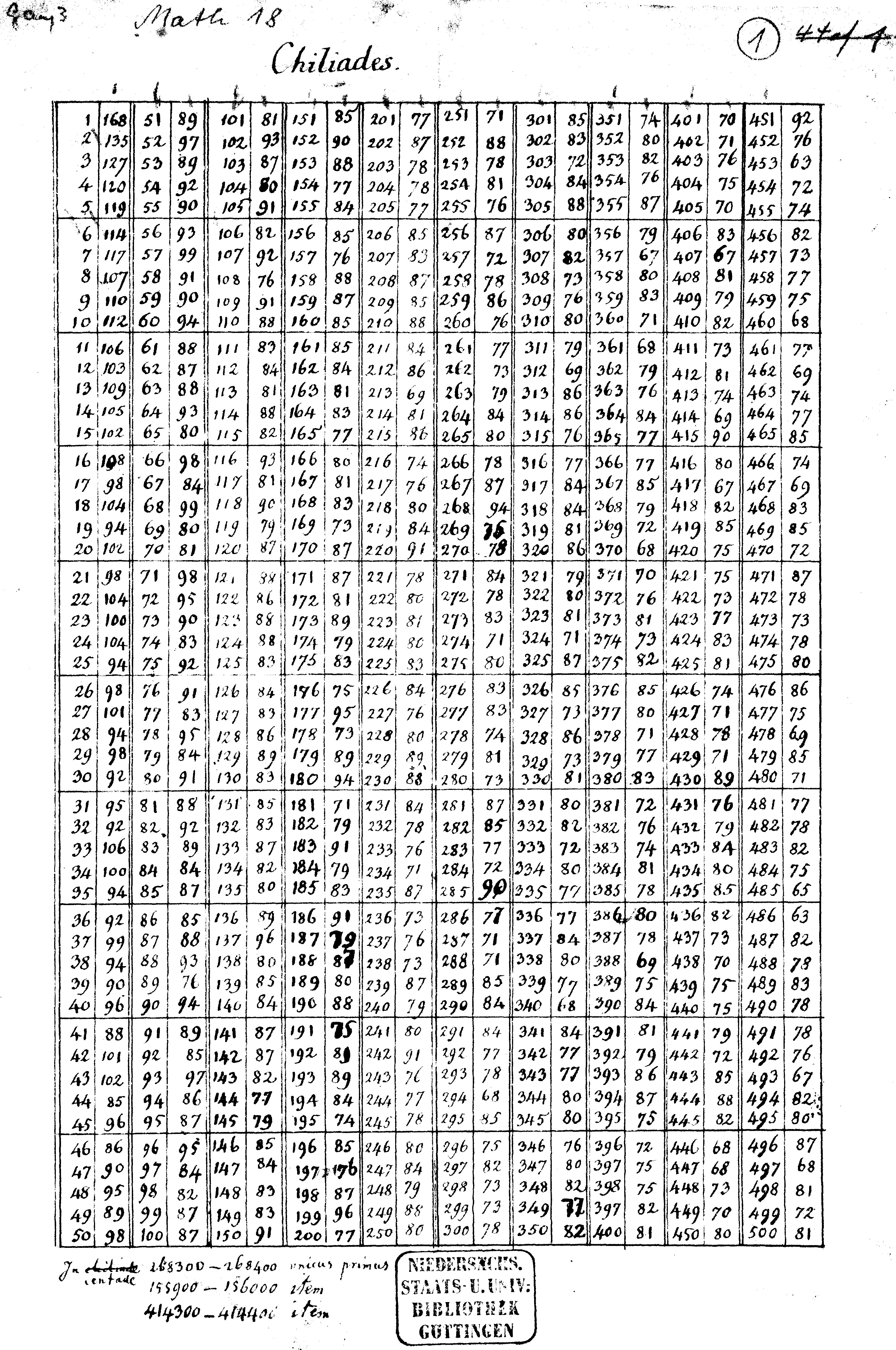}\quad\quad
	\end{subfigure}\quad\quad\quad
	\begin{subfigure}{0.45\linewidth}
		\includegraphics[scale=0.355]{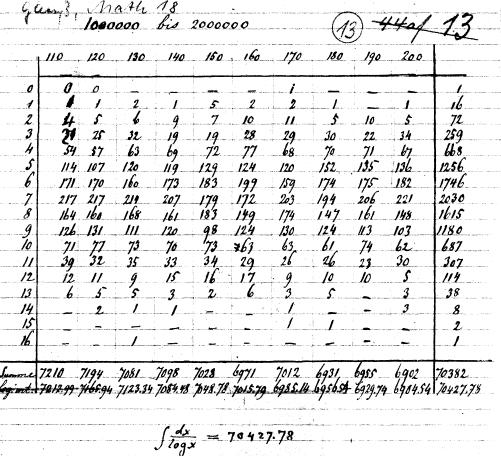}
	\end{subfigure}
	\caption{Hojas manuscritas por Gauss fechadas antes de 1840}
\end{figure}

Gauss, como hemos visto anteriormente, era capaz de estimar numéricamente la funcion $li(x)$ para valores muy grandes de $x$.  Para compararla con la función $\pi(x)$, tal y como Gauss deseaba, era obligado ampliar las tablas de primos. En la mencionada carta a Dase le pide a éste que amplíe de los 6 a los 10 millones. Dase cumplió el encargo publicándose de manera póstuma en 1862-65, \cite{Das1862}, y que Gauss tampoco llegó a ver.

\subsection*{B. Riemann y la función $li(x)$}

El 3 de noviembre de 1859, B. Riemann (1826-1866) comunicó a la Academia de Berlín un breve pero trascendental documento \textit{"Über die Anzahl der Primzahlen unter einer gegebenen Grösse"}, \cite{Rie1860}. 
Será el único trabajo que publicó sobre Teoría de Números.

\textit{No creo poder expresar mejor mi agradecimiento por la distinción que la Acade-
mia me ha hecho al nombrarme como uno de sus correspondientes que haciendo uso
del privilegio que conlleva dicho nombramiento para comunicarle una investigación
sobre la densidad de los números primos. Una materia que a causa del interés que
Gauss y Dirichlet le han dedicado durante muchos años no parece indigna de una
tal comunicación.} 

Riemann demuestra la famosa identidad que contecta $li(x)$ y $\pi(x)$,
$$\pi(x) +\frac{1}{2}\pi(x^{1/2}) + \frac{1}{3}\pi(x^{1/3}) + … = li(x)-\sum_{\rho}li(x^{\rho})+ E(x)$$
donde $E(x)$ es un término de error y $\rho$ recorre los ceros de la zeta de Riemann, dando lugar a la formulación de la celebrada Hipótesis de Riemann. Al final del artículo Riemann observa que \textit{... la comparación de $li(x)$ con el número de primos menores que $x$ realizada por Gauss y Goldschmidt hasta tres millones,
	muestra que el número de primos es menor que $li(x)$ en los primeros cientos de miles...}\smallskip
Y propone como consecuencia de su fórmula anterior la mejora 
$$\pi(x)\sim li(x)-\frac{1}{2}li(x^{1/2})-\frac{1}{3}li(x^{1/3})-\frac{1}{5}li(x^{1/5})+\frac{1}{6}li(x^{1/6})-\frac{1}{7}li(x^{1/7})...$$

\begin{figure}[htb]
	\centering
	\includegraphics[scale=0.225]{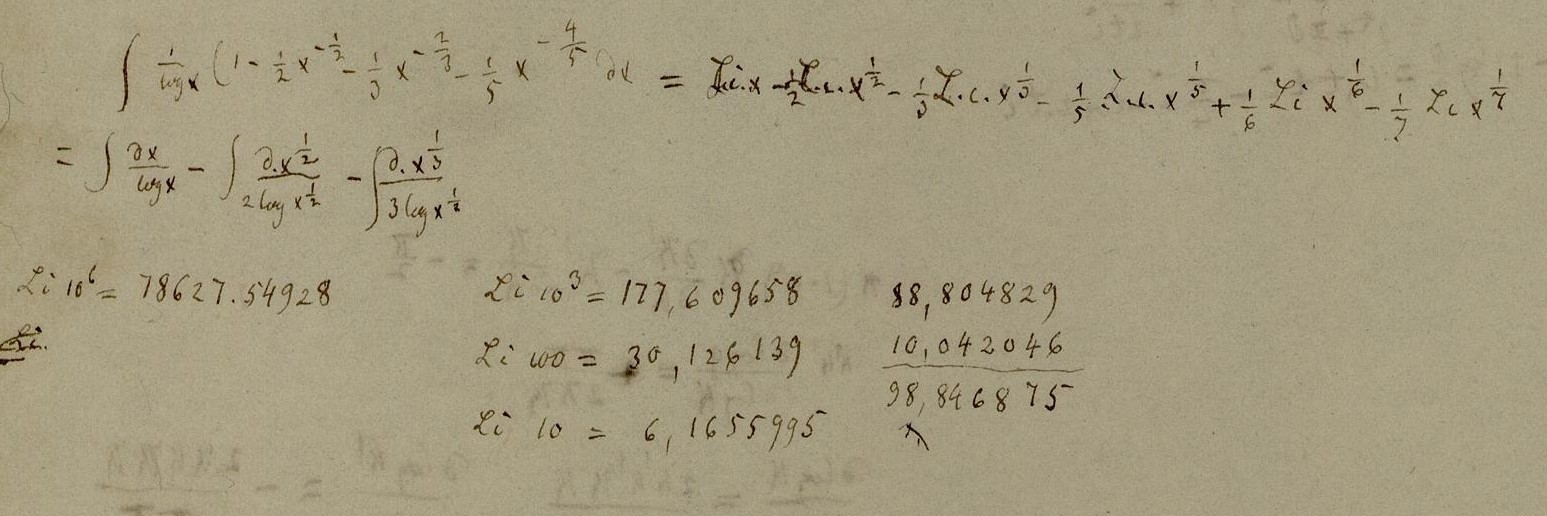}\medskip
	\caption{Estimacion de $Li(x)$ en los Nachlass de Riemann}
	\label{fig:Riemann}
\end{figure}

 Riemann realizó (en contra de lo que se ha pensado durante mucho tiempo)  cálculos \textit{numéricos} con la función $li(x)$. Como observamos en el manuscrito, \cite{Rie1800}, ver  Figura \ref{fig:Riemann}, Riemann estima numéricamente entre otras $\pi(10^6)$ mediante el cálculo explícito de los valores de $li(x)$ necesarios para probar la validez de su fórmula de aproximación.\footnote{Riemann denota el logaritmo integral con mayúscula, $Li(x)$. En la actualidad se utilizan ambas notaciones $li(x)$ y $Li(x)$ indistintamente, \cite{Der2003}. } \smallskip

\noindent Es esta la primera vez que aparece \textit{publicado en una revista científica} el nombre de Gauss (y de Goldschmidt) en relación al tema de la distribución de los números primos. Algunos autores, \cite{Pat2017}, destacan como posible germen de este trabajo la relación de Gauss, Dirichlet y Goldschmidt con Riemann; se sabe además que la relación entre estos tres últimos era estrecha.

\section{La carta que reescribe la historia.}

E. Landau en la introducci\'on hist\'orica de su libro \textit{“Handbuch der Lehre von der Verteilung der Primzahlen”} de 1909, \cite{Lan1909}, coloca a Gauss, quiz\'a de manera sorprendente a los ojos de un lector actual,  despu\'es de Riemann justific\'andolo de esta forma: \bigskip

\textit{“Yo hablo de Gauss s\'olo ahora porque \'el mismo no tiene nada publicado sobre la teor\'{\i}a de los n\'umeros primos; los detalles de su compromiso con el tema fueron descubiertos en 1863, despu\'es de la publicaci\'on del tratado de Riemann, en el Volumen II de la recopilacion de sus obras, \cite{Gau1863}, donde aparece una carta muy interesante a Encke que se fecha el 24 de diciembre de 1849….”}  \bigskip

La afirmación de Landau no es del todo correcta. La primera referencia al papel de Gauss es un poco anterior, de 1859. En el libro de “Exercices D’Analyse Numerique” de V. A. Le Besgue,\cite{Leb1859}, en el apartado dedicado a la \textit{"enumération des nombres premiers"} cita las aportaciones de Legendre, Dirichet y Tchebycheff. Pero también incluye las aportaciones de Gauss guiándose en la reciente aparición de la correspondencia entre Bessel y Olbers publicada en 1852. 

Es importante destacar que la famosa carta a la que alude Landau, es de hecho la respuesta a otra enviada por J.F.Encke (1791-1865) a Gauss el 4 de diciembre de 1849 \footnote{De esta carta se conocía tan solo el fragmento que se encuentra en notas finales de las Werke t. II,  pero no había sido publicada en su totalidad hasta 2018 en el epistolario entre Gauss y Encke.}, \cite{Enk2018}, \cite{Gauss}. La misma revela algún detalle interesante.

Encke comenta a Gauss que unos meses atrás ha asistido a una conferencia-seminario de Dirichlet en la Academia de Berlín en la que menciona \textit{la relación entre el logaritmo hiperbólico y el número de primos hasta un número dado, que Legendre había propuesto por primera vez... Esto me llevó a buscar alguna fórmula lo mas simple posible que se acercara más al numero de primos.} Encke propone la expresión $\frac{n}{\log n}10^{\frac{1}{2\log n}}$ a la que acompaña con una tabla.\medskip

En el contexto historico, nos detenemos a destacar la oportunidad de la carta de Encke compartiendo sus ideas acerca de la distribución de los numeros primos. Si ésta no se hubiese producido o Gauss simplemente no la hubiese contestado ya por falta de tiempo o interés en ese momento, los conocimientos que Gauss tenía en este área seguramente (puesto que no escribió nada más sobre la cuestión) nunca se habrían mostrado al mundo, perdiéndose algunas de sus contribuciones. 

La carta de Gauss a Encke del 24 de diciembre de 1849 se encuentra reproducida en el volumen II de los Werke, a continuación de la enumeración de los números primos hasta los tres millones realizada por Gauss y Goldschmidt\footnote{El editor Dr. Schering, afirma que los resultados para el primer millón estaban manuscritos por Gauss, y los del segundo y tercer millones estaban íntegramente escritos por Goldschmidt.}, \cite{Gau1863}. 
Su contenido cambiará el relato temporal de los acontecimientos tal y como describe con gran detalle el conocido artículo de L.J. Goldstein,  \cite{Gol1973}.

Gauss responde a Encke, \cite{Gauss}: 

 \textit{"Mi distinguido amigo: Tus observaciones sobre la frecuencia de los primos fueron de interés para mí de varias formas. Me has recordado mis investigaciones sobre el tema, que se remontan a un pasado muy lejano en 1792 o 1793, después de que adquiriera los suplementos a las tablas logarítmicas de Lambert. Incluso antes de haber comenzado mis investigaciones mas detalladas en aritmética superior, uno de mis primeros proyectos fue centrar mi atención en la frecuencia decreciente de los primos, para lo cual había enumerado los primos en varias chiliadas \footnote{"De chiliada en chiliada" significa "de 1000 en 1000"} y apuntado los resultados en páginas en blanco adjuntas.}

\textit{"Pronto reconocí que, detrás de todas sus fluctuaciones, la frecuencia es en promedio inversamente proporcional al logaritmo, de modo que el número de primos inferiores a una cantidad dada $n$ es aproximadamente igual a 
$$\int\frac{dn}{log n}$$
donde el logaritmo se entiende hiperbólico. Posteriormente, cuando conocí las tablas de Vega\footnote{La referencia que Gauss hace a la tabla de Vega se publica en realidad en febrero de 1797.}  (1796), que tenían una lista hasta 400.031, amplié mis cálculos, confirmando esa estimación. En 1811, la aparición del Cribrum de Chernac me dió gran satisfacción y con a menudo (pues me faltaba paciencia para hacerlo continuamente) empleaba cuartos de hora desocupados para contar algunas chiliadas; aunque finalmente lo dejé de lado sin completar un millón. No fue hasta más adelante que hice uso de la diligencia de Goldschmidt para llenar los vacíos que quedaron en el primer millón y para continuar el cálculo de acuerdo con las tablas de Burckhardt. Así (ya hace años de esto) 
la enumeración se completó para los primeros tres millones y se contrastó con la integral. Un breve extracto es el siguiente, \cite{Gauss}}

\begin{figure}[htb]
	\centering
	\quad\quad\includegraphics[scale=0.1]{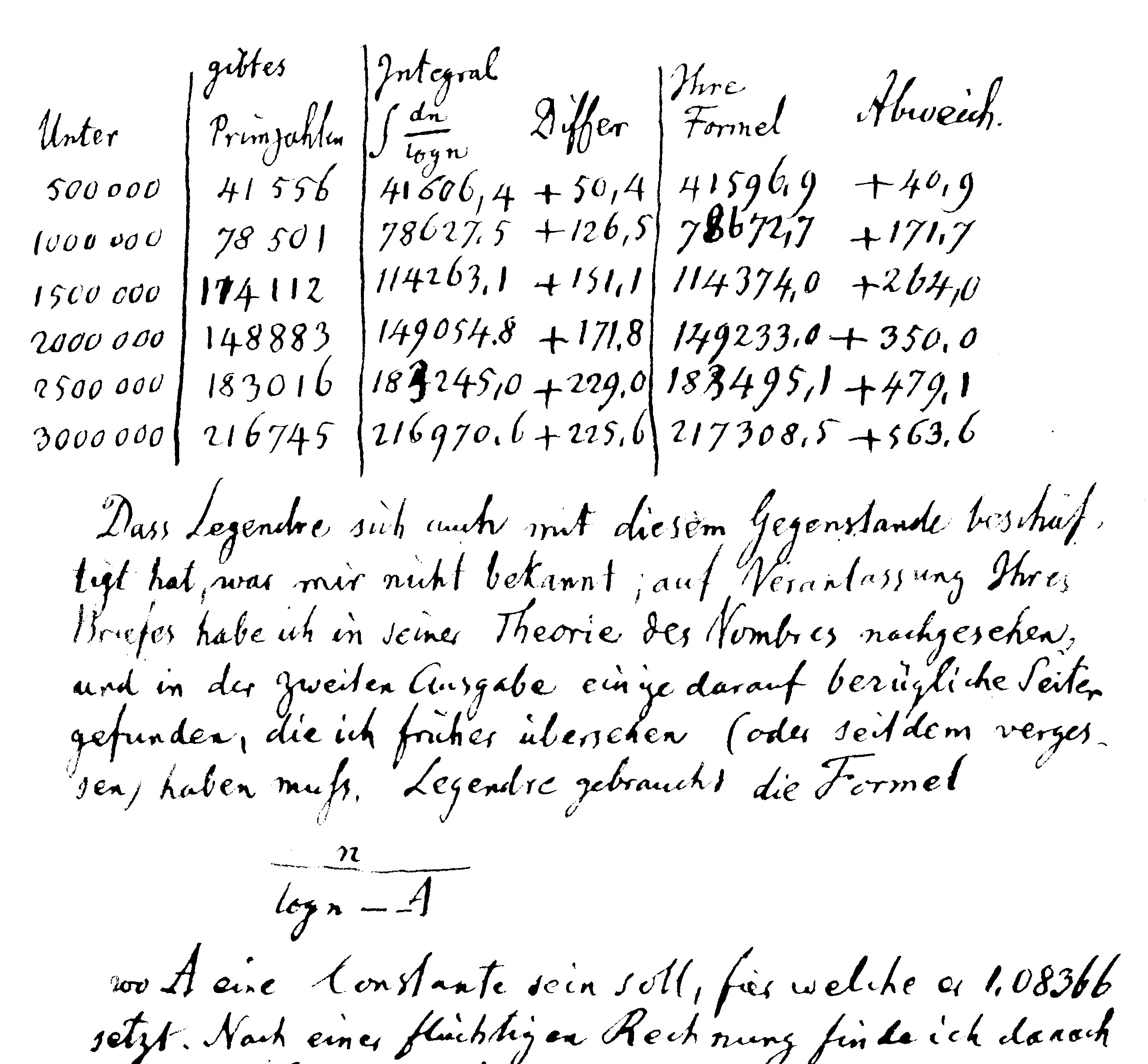}
	\caption{Extracto carta de Gauss, 1849}
\end{figure}


\textit{No era consciente de que Legendre también había trabajado en este tema; tu carta me hizo mirar en su Théorie des Nombres, y en la segunda edición encontré unas pocas páginas que debí haber pasado por alto (o ahora olvidado). Legendre usaba la fórmula $$\frac{n}{\log n -A}$$ donde $A$ es la constante que toma valor 1'08366.}

 Gauss compara los datos de Legendre con sus propios del logaritmo integral, observando que para los rangos de $n$ hasta 3.000.000, los de Legendre se ajustan mejor. Tantea otros posibles valores de $A$ y \textit{no me atrevo a conjeturar que si el valor límite \text{(de $A$)} cuando n tiende a infinito es 1 o distinto de 1.} \footnote{Como hemos visto anteriormente, Tchebycheff acababa de demostrar que $A=1$.}

Gauss finaliza la carta advirtiendo de algunas discrepancias entre la enumeración hecha tanto por Encke como por Lambert y su propios resultados. Alude también a Dase (en la que probablemente sea la primera mención de su nombre en conexión con tablas de factores) y anima a Encke para que hable con éste y continúe calculando más tablas de primos. Este intercambio continuará en sucesivas cartas (21 de febrero y 11 de marzo de 1851) \cite{Enk2018}.

\section{Discusión y evidencias}

Lo expuesto hasta ahora evidencia que hay un desfase entre lo que se publica y cuándo se escribe. Por un lado la comunidad matemática tiene noticias de los trabajos que se van publicando en una secuencia temporal dada. Por otro tenemos lo que podemos denominar el “polymath alemán” liderado por Gauss, Dirichlet... que maneja, además de lo publicado, el conocimiento que emerge en su entorno académico privado.

Es conocido que la forma de trabajar de Gauss \textit{le llevaba a dar a sus investigaciones la forma de perfectas obras de arte...y nunca dió un trabajo para publicar hasta darle la perfección que deseaba...} siendo su  lema\footnote{La frase aparece, junto con numerosas noticias personales, en el libro-memorial de Sartorius von Waltershausen, profesor en Götinga y amigo de Gauss, escrito en 1856 a la muerte de éste.  La traducción al inglés de "Carl Fiedrich Gauss, a memorial" fue realizada por Helen Worthington Gauss, biznieta del matemático, y publicada en 1966.
	Ha sido reimpresa recientemente (2018) por www.forgottenbooks.com}: \textit{"Pauca, sed matura"} (Pocos, pero maduros). \cite{Wal1856}

Landau, \cite{Lan1909}, se muestra crítico al valorar el relato que escribe Gauss en su carta a Encke y dice de élla que  \textit{...pero Gauss no ha demostrado nada; él tiene esa carta escrita en la vejez con 72 años, hacia el final de su vida y sólo contiene una explicación y extensión de las presunciones que Gauss había tenido de niño}. 

Revisando la carta de Gauss a Encke podemos ahora contrastar algunas de las afirmaciones incluídas en la misma así como ciertas llamativas omisiones:

\textbf{1. } Gauss dice haber pasado por alto u olvidado la segunda edición del libro de Legendre. Sin embargo el mismo Gauss en el prefacio de sus “Disquisitiones Arithmeticae” (1801), \cite{Gau1801}, alaba la primera edición del sobresaliente trabajo de Legendre (1797-98). Además hoy sabemos que Legendre le mandó una copia (de la segunda edición) junto con una carta fechada el 8 de noviembre de 1808, \cite{Gauss}, donde le señala el capítulo 8 de la parte cuarta del libro que es precisamente donde se encuentra publicada su famosa fórmula de distribución de primos así como la tabla de primos que hemos expuesto anteriormente (ver figura 1).  No hay constancia en \cite{Gauss}, que contestase. Esta edición era, como hemos explicado, referencia para toda la comunidad matematica del momento. Para ilustrar su popularidad recordemos el comentario de N. Abel en 1823 a su profesor Holmboe, refiriéndose a la fórmula de Legendre \textit{que es probablemente la fórmula más importante en todas las matemáticas}. Asimismo, como hemos indicado anteriormente,  Dirichlet en 1838 envía a Gauss una nota manuscrita sobre tal fórmula.

\textbf{2. } Después de lo expuesto en secciones anteriores sobre la correspondencia mantenida entre Gauss y Bessel, sorprende la nula mención a éste último. Es conocido que la relacion epistolar entre ambos comienza en 1804 por mediacion de Olbers, \cite{Gau2005}. 
También es sabido, y podría explicar la omisión a Bessel, que el tono amable de sus cartas cambia a partir de la decada de los 30 a raíz de una discusión que ambos mantuvieron en la que Bessel reprocha a Gauss el no comunicar sus progresos matematicos \textit{"nunca reconociste la obligación de comunicar parte tus investigaciones para promover el conocimiento real de los temas y así ellos viven en la noche del mundo"}. Esta falta de respeto que Gauss consideró inaceptable, como le comentó a Schumacher \cite{Sch1860}, le llevó a interrumpir su correspondencia durante cinco años y medio (entre 1834 y 1839), mientras que Bessel le seguía enviando regularmente alguna carta con material relativo a los trabajos que iba realizando. Al morir Bessel en 1846, con la correspondencia Bessel-Gauss preparada publicar (196 cartas), Gauss toma la decisión de no incluir la carta de la discusion con Bessel, lo que detiene el proyecto. Posteriormente, con los dos ya fallecidos se acabará por incluir en 1880, \cite{Eng1880}, \cite{Gau2005}.

\textbf{3. } Gauss recuerda sus enumeraciones de primos en chiliadas en su juventud. Más aún, como le comunica por carta a Bessel en 1811 \cite{Gauss}, calcula en dicho año los primos por miriadas y además sabemos que el rango en el que trabaja es hasta el primer millón (sin llegar a completarlo). La referencia que hace a Goldschmidt, la podemos documentar hacia en la decada de los años 30 y principios de los 40. C. W. B. Goldschmidt (1807-1851) fue ayudante de Gauss en su observatorio desde 19 el diciembre 1834 y profesor de B. Riemann. Los datos recopilados por Gauss y Goldschmidt sobre el crecimiento de la integral logarítmica en comparación con la $\pi(x)$ fueron citados en el célebre artículo de Riemann, \cite{Rie1860}. 

\textbf{4. } Curiosamente, en el año 1849 confluyen algunos acontecimientos que creemos que pueden enmarcar la famosa respuesta de Gauss a Encke. Por un lado la presentación de la tesis de Tchebycheff, \cite{Tch1848}, el 24 de mayo de 1848 y por otro la aparición en revista científica del artículo de Hargreave \cite{Har1849} en el verano de 1849. A este respecto, S.J. Patterson de la Universidad de Götinga, \cite{Pat2017}, escribe \textit{"no sabemos si Gauss conocía el trabajo de Tchebycheff pero es conocido que Gauss tenía numerosos contactos con la Academia de San Petersburgo"}. En esa época además era capaz de leer y hablar ruso con facilidad \cite{Wal1856}.  Patterson confirma también que \textit{Philosophical Magazine}, la revista en la que había aparecido el artículo de Hargreave, \textit{"llegaba asiduamente a Götinga y que al menos Riemann lo leía. Podemos concluir ... que Gauss tuvo oportunidad de verlo. Por otro lado, Tchebycheff visitará Berlin en el verano de 1852 invitado por Dirichlet ... podemos esperar que Riemann tuviera conocimiento de estos avances"}. 

\textbf{5. } Por la carta que Gauss envía a Dase en 1850, \cite{Das1862}, el interés de aquel por que se ampliasen las tablas de primos coincide en el tiempo con la observación que hace Tchebycheff, \cite{Tch1848}, de que hasta que no se llegue a tener datos por encima de los 10 millones no se podrán observar grandes diferencias entre la fórmula propuesta por Legendre y la de Tchebycheff-Gauss. Con el tiempo se comprobará que de hecho será más precisa esta última a partir de los 5 millones.

\textbf{6. } La frase inicial de Gauss  \textit{mis investigaciones sobre el tema...se remontan a un pasado muy lejano en 1792 o 1793} ha dado lugar a diversas interpretaciones, provocando incluso que se haya llegado a adjudicar a Gauss, \cite{Gol1973},  un conocimiento de la relación entre $li(x)$ y $pi(x)$ demasiado elaborado para una época tan temprana.

Hoy sabemos que las numerosas leyes asintóticas que Gauss había encontrado a través de una combinación de consideraciones heurísticas y teóricas fueron publicadas póstumamente como lista de proposiciones, \cite{Buh1981} (ver \cite{Dun2004}). Concretamente en 1917 en el tomo 10, volumen 1 de los Werke, \cite{Gau1917} se recogen apuntes manuscritos que Gauss anotaba en las hojas en blanco de los libros que iba adquiriendo. 

En la última pagina del ejemplar de tablas de logaritmos de Schulze, \cite{Schu1778} que Gauss poseía desde los 14 años, aparece entre otras anotaciones escrito en alemán y fechado en mayo de 1796:

\centerline{\textit{Primzahlen unter a(=infinito) $\frac{a}{la}$}}

seguido de 

\centerline{\textit{Zahlen aus zwei Factoren $\frac{lla \cdot a}{la}$}} 

Asimismo, el mismo volumen de los Werke incluye el famoso diario de Gauss, que el denominó Notizenjournal, que consta de 19 páginas y escribió entre los años 1796 y 1814 donde, como él mismo reconoce, \textit{"apuntaba la gran cantidad de ideas que le venían a la mente, para que no se le olvidasen y poder volver sobre ellas"}. En total contiene 146 anotaciones breves, ocasionalmente difíciles, si no imposibles, de interpretar correctamente. En la anotación 9 fechada el 31 de mayo de 1796 se lee: {\textit{"Comparationes infinitorum in numeris primis et factoribus contentorum."}}

Estas dos notas por la coincidencia de las fechas son tomadas como juntas por numerosos autores entre ellos P. Bachmann,\cite{Gau1917}, que fue uno de los encargados por los editores de escribir sobre los trabajos de Gauss en el volumen 2 del Werke tomo 10 y que como él indica: \textit{Gauss anticipó alguna de las leyes asintóticas de los números primos, continuó trabajando de vez en cuando en ello después de 1811 y luego dejó que Goldschmidt continuara extendiendo las tablas que aparecen en el tomo 2 de los Werke; además pronto tuvo conocimiento del logaritmo integral como la funcion asintótica adecuada para ser la expresión para el conjunto de números primos bajo un límite dado.} Es claro que Gauss sabia muchas cosas que lamentablemente no plasmo con la suficiente claridad, pero como indica Bachmann este conocimiento lo fue adquiriendo a lo largo del tiempo.

Por otro lado, en los años 1792-93 no existían datos (tablas) suficientes. Las tablas de Schulze de 1778 que Gauss posee, únicamente contienen los primos hasta 10009.  Las tablas de primos existentes no superan a finales del XVIII los 400.000, llegando al millón en el trabajo de Chernac de 1811.
Por otro lado el cálculo de la función $li(x)$ para valores grandes hasta $x=1000000$ llegará con Bessel en 1810.

Para concluir, y a modo de resumen ilustrativo presentamos la siguiente tabla, inspirada en \cite{Jam2003},  que muestra los datos numéricos aportados por cada una de las expresiones que se han discutido a lo largo de todo el artículo. Tanto fórmulas como cálculos han sido protagonistas imprescindibles en el difícil acercamiento al problema de la distribución de los números primos. 

\bigskip

\begin{center}
	\begin{tabular}{| c | c | c | c | c |}
		\hline
		\multicolumn{5}{ |c| }{Tabla comparativa.  } \\ \hline
		n & $\pi(n)$ & $\frac{n}{\log n}$ & $\frac{n}{\log n -1,08366}$  & $\int_2^n\frac{dx}{\log x}$  \\ \hline
		1000 & 168 & 145 &  172 &  177  \\
		10000 & 1229 & 1086 & 1230 & 1246  \\
		50000 & 5133 & 4621 & 5136 & 5166  \\
		100000 & 9592 & 8686 & 9588 & 9630  \\
		500000 & 41538 & 38103 & 41533 & 41607  \\ 
		1000000 & 78498 & 72382 & 78543 & 78628  \\ 
		10000000 & 664579 & 620421 & 665140  & 664918  \\ \hline
	\end{tabular}
	\label{tab:Jameson}
\end{center}
\bigskip

\textbf{Comentario a la bibliografía. } Queremos destacar de manera especial dos referencias de obligada consulta y que nos sirvieron como punto de partida en nuestra investigación. La primera es el “Handbuch der Lehre von der Verteilung der Primzahlen” de Edmund Landau (1877-1938) en 1909, \cite{Lan1909}, posiblemente ha sido la que más ha influido en los especialistas en números primos del siglo pasado siendo alabada por figuras tan destacadas en teoría de numeros como Hardy, Littlewood y Montgomery... Es la primera obra que recoge todo lo conocido hasta la fecha de publicación acerca de la distribución de los números primos. Es extremadamente detallado, presenta de forma sistemática, clara y atractiva la materia, aportando valiosísimas referencias para poder seguir las líneas de investigación de los distintos autores, si tiene algún inconveniente es que está en alemán.

La siguiente referencia es la de W. Narkiewicz (1936- ) que publica en el año 2000 “The Development of Prime Number Theory: from Euclid to Hardy and Littlewood”, \cite{Nar2000}. Este texto recoge desde la antigüedad hasta la primera decada del siglo XX, el desarrollo de la teoría de la distribución de los números primos. Presenta los principales resultados con sus demostraciones y también da, en multitud de pequeños comentarios, una visión general del desarrollo en los ultimos 80 años. Contiene más de 1800 referencias. 

Es justo también citar la guía proporcionada por \cite{Har1938}, \cite{Edw1974} o \cite{Der2003}.

\end{document}